\documentclass[12pt]{article} 

\usepackage{amsmath,amsthm,amssymb,latexsym}
\usepackage{amsfonts,mathrsfs,eucal}
\usepackage{enumerate,subfigure,verbatim}
\usepackage{graphics,graphicx,tabularx,cite}

\usepackage{tikz}
\usetikzlibrary{arrows}
\usepackage[all,cmtip,2cell]{xy} \UseAllTwocells 

\usepackage[margin=1in]{geometry}
\usepackage[latin1]{inputenc}

\theoremstyle{definition}
\newtheorem{thm}{Theorem}[section]

\newtheorem{pro}[thm]{Proposition} 
\newtheorem{lem}[thm]{Lemma}

\newtheorem{cor}[thm]{Corollary}

\newtheorem{AshleyExt}[thm]{Dold-Kan Extension}

\newtheorem*{mainthm*}{Theorem 7.1}

\newtheorem{Def}[thm]{Definition}
\theoremstyle{remark}
\newtheorem{rmk}[thm]{Remark}
\newtheorem{ex}[thm]{Example}

\newcommand{\sgrp}{\textbf{sGrp}}
\newcommand{\sgrpd}{\textbf{sGrpd}}

\newcommand{\Top}{\textbf{Top}}

\newcommand{\xm}{\textbf{xm}}

\newcommand{\Xc}{\textbf{Xc}}
\newcommand{\XC}{\textbf{XC}}
\newcommand{\xc}{\textbf{xc}}

\newcommand{\rhom}{\textbf{Rhom}}

\newcommand{\crs}{\mathsf}
\newcommand{\one}{\crs{1}}

\newcommand{\B}{\crs{B}}
\newcommand{\C}{\crs{C}}
\newcommand{\D}{\crs{D}}
\newcommand{\E}{\crs{E}}
\newcommand{\G}{\crs{G}}
\renewcommand{\H}{\crs{H}}
\newcommand{\I}{\crs{I}}

\newcommand{\Q}{\crs{Q}}
\newcommand{\N}{\crs{N}}

\newcommand{\cat}{\mathcal}

\newcommand{\cG}{\cat{G}}

\renewcommand{\1}{\textbf{1}}
\renewcommand{\i}{^{-1}}
\newcommand{\me}{^{\prime}}

\newcommand{\coker}{\text{coker~}}
\renewcommand{\Im}{\text{Im~}}

\newcommand{\p}{\partial}
\renewcommand{\d}{\delta}

\newcommand{\ol}{\overline}
\renewcommand{\bf}{\textbf}
\renewcommand{\it}{\textit}
\renewcommand{\t}{\text}

\begin{document}

\title{n-Butterflies: Algebraically Modeling Morphisms between Homotopy $n$-Types}
\date{}
\author{ Ivan Dungan \\ Department of Mathematical Sciences \\ U.S. Military Academy, West Point}  

\maketitle

\abstract{Crossed modules are known to be a model of pointed connected homotopy $2$-types; formally, the homotopy category of crossed modules is equivalent to the category of pointed connected homotopy 2-types.  In forming the homotopy category of crossed modules, one must resort to computing derived morphisms using non-constructive topological methods, but Behrang Noohi was able to find an algebraic model of these derived morphisms called butterflies.  The result is a completely algebraic model of pointed connected homotopy 2-types.   Reduced crossed complexes are a generalization of crossed modules and model a subclass of pointed connected homotopy types.  We will present algebraic objects called n-butterflies which satisfy similar properties to butterflies and begin to generalize the theory of butterflies to model morphisms of pointed connected homotopy n-types.}

\tableofcontents

\section{Introduction}

A \it{crossed module} $[\C:\p]$ is a group homomorphism $\p:\C_2\rightarrow \C_1$ with a right action of $\C_1$ on $\C_2$, denoted by $a^x$, which satisfies the axioms:
\begin{align*}
&\text{CM1}  ~~~ \p(a^x) =x^{-1}\p(a)x  \text{ for all } x\in \C_1, a\in \C_2 \\ 
&\text{CM2}  ~~~ a^{\p(b)} =b^{-1}ab \text{ for all } a,b\in \C_2 
\end{align*}
These imply that the action descends to conjugation under $\p$ and lifts to conjugation when restricted to $\p(\C_2)$, respectively.

Crossed modules form a model category and the resulting homotopy category is equivalent to the category of pointed connected homotopy $2$-types.  Further, the morphisms of the homotopy category can be described algebraically by considering objects called butterflies as morphisms \cite{Noohi2}.  Specifically, a \it{butterly} $B=(\E,p,f,\alpha, \beta )$ from a crossed module $[\H:\p]$ to another $[\G:\d]$ is a commutative diagram
\begin{equation}
\xymatrix@=1em{ \H_2 \ar[dr]^\alpha \ar[dd]_\p & &\G_2 \ar[dd]^\d \ar[dl]_\beta \\ &\E \ar[dr]_f \ar[dl]^p  \\ \H_1  && \G_1  }
\end{equation}
where both diagonal sequences are complexes and $\G_2\rightarrow \E \rightarrow \H_1$ is short exact.

It was shown that butterflies model morphisms of the homotopy category of crossed modules. The resulting category of crossed modules and butterflies as morphisms is a completely algebraic model of pointed connected homotopy $2$-types.  This work was extended to over a site by the first author and Noohi in \cite{AldrovandiBehrang1}.

Similar to crossed modules, crossed complexes admit a model category structure whose homotopy category models a subclass of homotopy types \cite{BrownGolasinski1}.  The aim of this paper is to present algebraic objects called $n$-butterflies which satisfy similar properties to butterflies and begin to generalize the theory of butterflies in order to have a completely algebraic model.  We will focus our attention on reduced $n$-crossed complexes which model a subclass of pointed connected homotopy $n$-types and begin with the theory of these objects.  

We follow the procedure of computing butterflies by first defining a $n$-pushout of a morphism of crossed complexes.  The $n$-pushout applied to a derived morphism of crossed complexes produces a new fraction which can be expressed as a diagram whose shape dignifies its name, $n$-butterfly.  We will define special $n$-butterflies which form a groupoid denoted by $n\B_\Q(\H,\G)$ where $\H,\G$ are reduced $n$-crossed complexes.  These objects are less topological and more algebraic than derived morphisms, however, their algebraic structure encodes enough information to still characterize the morphisms of the homotopy category.  We have the following main result.

\begin{mainthm*}
Let $[\H:\p]$ and $[\G:\d]$ be reduced $n$-crossed complexes.  Then we have the bijection $[\H,\G]_{\Xc} \cong \pi_0(n\B_\Q(\H,\G))$.
\end{mainthm*}

The following theory appears to be adaptable to crossed complexes; however, we will postpone this excursion to a later paper and only concentrate on reduced crossed complexes.  We will now review some necessary background on the homotopy theory of crossed complexes which can be skipped if already known.

\section{Crossed Complexes and Homotopy Types}

Below we give a quick summary of the results needed for our purpose beginning with the definition of crossed complexes.  For further details, see \cite{BrownHigginsSivera1} for an initial resource.

\begin{Def}
A \it{crossed complex} $[\C:\p]$ is a chain of groupoids over a set of objects $\C_0$ represented by the following diagram and satisfying the axioms below.
\begin{equation} \label{xcdiagram}
\xymatrix{ \ldots \ar[r] &\C_{n+1}\ar[r]^{\p}  &\C_n \ar[r]^{\p} &\ldots \ar[r]^{\p} &\C_2 \ar[r]^{\p} &\C_1 \ar@<.5ex>[r]^{\p_1} \ar@<-.5ex>[r]_{\p_0} &\C_0 }
\end{equation}
\begin{enumerate}
\item [\bf{Xc1}] For $k\geq 2$, $\C_k$ is a discrete groupoid over $\C_0$ where $\C_k(x)$ denotes the group of maps $\C_k(x,x)$.
\item [\bf{Xc2}] For $k\geq 3$, the groups of $\C_k$ are abelian.
\item [\bf{Xc3}] For $x,y\in \C_0$, $\C_1(x,y)$ acts on the right of $\C_k(x)$ which is conjugation when $k=1$.
\item [\bf{Xc4}] The maps $\p_0,\p_1$ are the source and target maps of $\C_1$, and $\p$ is a chain map of groups over each object in $\C_0$ which preserves the action of $\C_1(x,y)$.
\item [\bf{Xc5}] For $x\in \C_0$ and $u\in \C_2(x)$, $\p(u)$ acts trivially on $\C_k(x)$ for $k\geq 3$ and by conjugation on $\C_2(x)$.
\end{enumerate}
A morphism of crossed complexes is a family of morphism of groupoids $\{f_k:\C_k\rightarrow \D_k\}$ for $k\geq 1$ which respect a morphism of objects $f_0:\C_0\rightarrow \D_0$ such that the diagram
\begin{equation} \label{xcmorphism}
\xymatrix{ \ldots \ar[r] &\C_{n+1}\ar[r]^{\p} \ar[d]_{f_{n+1}} &\C_n \ar[r]^{\p} \ar[d]_{f_n}&\ldots \ar[r]^{\p} &\C_2 \ar[r]^{\p} \ar[d]_{f_2}&\C_1 \ar@<.5ex>[r] \ar@<-.5ex>[r] \ar[d]_{f_1}&\C_0   \ar[d]_{f_0}
\\  \ldots \ar[r] &\D_{n+1}\ar[r]^{\d}  &\D_n \ar[r]^{\d} &\ldots \ar[r]^{\d} &\D_2 \ar[r]^{\d} &\D_1 \ar@<.5ex>[r] \ar@<-.5ex>[r] &\D_0 }
\end{equation}
commutes and is compatible with the action of $\C_1$ on $\C_k$.  Crossed complexes and these morphisms form a category denoted by $\Xc$.  

\end{Def}

A \it{$n$-crossed complex} is a chain of groupoids of length $n$ satisfying the same properties as crossed complexes.  These objects and similarly defined morphisms form a category denoted by $n\Xc$ which is fully embeddable in $\Xc$.  Crossed complexes and $n$-crossed complexes with one object are called \it{reduced crossed complexes} and \it{reduced $n$-crossed complexes}, respectively, and form full subcategories $\xc$ and $n\xc$ of $\Xc$ and $n\Xc$, respectively. 

For reduced crossed complexes, we will usually ignore the object $*$ and visualize them as chains:
\begin{equation} \label{rxcdiagram}
\xymatrix{ \ldots \ar[r] &\C_{n+1}\ar[r]^{\p}  &\C_n \ar[r]^{\p} &\ldots \ar[r]^{\p} &\C_2 \ar[r]^{\p} &\C_1 }
\end{equation}
where the first homomorphism is a crossed module by axioms \Xc \bf{3}, \Xc \bf{4}, and \Xc \bf{5}.  We denote the category of crossed modules by $\xm$ and the diagram of the natural fully faithful functors $\xm \rightarrow \xc \rightarrow\Xc$.

These functors are examples of the $n$-skeleton functor, $sk^n$.   Coskeleton, truncation, and cotruncation functors can also be defined and are denoted by  $cosk^n$, $tr_n$, and $cotr_n$, resp.  

\begin{thm}\cite{BrownHigginsSivera1} \label{truncationadjunction}
The functors $sk^n,cosk^n,tr_n,cotr_n$ satisfy the adjoint\footnote{For the adjoint diagrams, left adjoints are written on the left and right adjoints are written on the right.} diagrams 
\begin{equation*}
sk^n:n\Xc \rightleftarrows \Xc:tr_n; \hspace{0.2in} tr_n:\Xc \rightleftarrows n\Xc:cosk^n; \hspace{0.2in} cotr_n:\Xc \rightleftarrows n\Xc:sk^n 
\end{equation*}
\end{thm}

In \cite{BrownHiggins2}, $n$-fold left homotopies are defined and are shown to form an internal path space $\XC(\C,\D)$.  Also, a tensor product is defined which together make the category of crossed complexes into a closed symmetric monoidal category $(\Xc, \otimes, \one)$ where $\one$ is the trivial reduced crossed complex called the \it{unital crossed complex}.  

\begin{Def}
Let $f,g:[\C:\p]\rightarrow [\D:\d]$ be two morphisms in $\Xc$.  A \it{homotopy} from $f$ to $g$ is a morphism of crossed complexes $h:\C\otimes \I \rightarrow \D$ which makes the diagram
\begin{equation*}
\xymatrix{ \C \ar[d]_{i_0} \ar[dr]^f \\ \C\otimes \I \ar[r]^h &\D \\ \C \ar[u]^{i_1} \ar[ur]_g }
\end{equation*}
commute.  If $h$ is a homotopy, then $f$ is \it{homotopic} to $g$ and denoted by $h:f\simeq g$.
\end{Def}

The following fact shows that homotopies as defined above are in fact equivalent to $1$-fold left homotopies.

\begin{thm}\cite{Tonk1} \label{eqofhomotopytheories}
Let $f,g:[\C:\p]\rightarrow [\D:\d]$ be two morphisms in \Xc.  A homotopy $h:f\simeq g$ is equivalent to defining a $1$-fold left homotopy $(g,\phi_n:\C_n\rightarrow \D_{n+1})$.  Moreover, $f$ is completely determined by the algebraic equations
\begin{align}
f_0(c_0)&=s(\phi_0(c_0))  \\
f_1(c_1)&=\phi_0(s(c_1))g_1(c_1)\d_{2}(\phi_1(c_1))\phi_0(t(c_1))\i  \\
f_n(c_n)^{\phi_0(t(c_n))}&=g_n(c_n)\d_{n+1}(\phi_n(c_n))\phi_{n-1}(\p_n(c_n))  &\text{for}~n\geq 2
\end{align}
\end{thm}

\begin{Def}
Let $f,g:[\C:\p]\rightarrow [\D:\d]$ be morphisms in $\xc$.  A \it{pointed homotopy} $h:f\simeq g$ is a homotopy which can be described by a pointed $1$-fold left homotopy.
\end{Def}

Explicitly, a homotopy $h:f\simeq g$ is pointed if and only if there exits a pointed left-homotopy $(g,\phi)$ in which $f$ is determined by the following.  % where \label{injhomotopy1} holds $f,g$ are morphisms in $n\xc$.
\begin{align} \label{pointeddetermine}
f_1(c_1)&=g_1(c_1)\d_2(\phi_1(c_1))  \\
f_k(c_k)&=g_k(c_k)\d_{k+1}(\phi_k(c_k))\phi_{k-1}(\p_k(c_k)) &\text{for~} k>1 
%f_n(c_n)&=\phi_{n-1}(\p_n(c_n)) \label{injhomotopy2}
\end{align}

\begin{Def}
Let $[\C:\p]$ be a crossed complex.  The \it{connected components} of $\C$ is the set $\pi_0(\C_1)$ of connected components of the groupoid $\C_1$.  For each $x\in \C_0$, the \it{fundamental homotopy group} of $\C$ at $x$ is the  group $\pi_1(\C,x)=coker[\d_2:\C(x)_2\rightarrow \C(x)_1]$.  The \it{$n$th-homotopy group} of $\C$ at $x$ is the group $\pi_n(\C,x)=H_n(\C(x))$.
\end{Def}

\begin{thm}
There exists a model structure on $\Xc$ where a morphism $f:[\C:\p]\rightarrow [\D:\d]$ in $\Xc$ is a \it{weak equivalence} if  the induced set map of connected components is a bijection and the induced maps of the fundamental homotopy group and $n$-homotopy groups are isomorphisms.  The map is a \it{fibration} if $f_1:\C_1\rightarrow \D_1$ is a fibration of groupoids\footnote{There is a detailed exposition on fibrations of groupoids in \cite{Brown2}.} and $f_n:\C_n(x)\rightarrow \D_n(f_1(x))$ is surjective for $n\geq 2$ and all $x\in \C_0$.
\end{thm}

\it{Cofibrations} for the model structure above are precisely the morphisms which have the left lifting property with respect to trivial fibrations (see \cite{GoerssJardine1}).  Since the unique map $\C\rightarrow \one$ where $\one$ is the trivial reduced crossed complex is clearly a fibration of crossed complexes for any $\C\in \Xc$, every crossed complex is fibrant; however, not every crossed complex is cofibrant.  A nice characterization of trivial fibrations is the following.

\begin{lem} \cite{BrownGolasinski1} \label{BGfibration}
A morphism $f:[\C:\p]\rightarrow [\D:\d]$ of crossed complexes is a trivial fibration if and only if the set function $f_0: \C_0 \rightarrow \D_0$ is surjective; the functor $f_1: \C_1 \rightarrow \D_1$ is full; for every $x\in \C_0$, the induced maps $\C_n(x)\rightarrow \ker \p_{n-1} (x) \times_{\ker \d_{n-1}(x) }\D_n(x)$ are surjective for $n\geq 2$.
\end{lem}

\begin{rmk} \label{cosknequalstwo}
For $n=2$ in the above theorem, $\ker \p_{n-1}(x)$ and $\ker \d_{n-1}(x)$ are the equalizers of the diagrams $\xymatrix{ \C_1(x)\ar@<.5ex>[r] \ar@<-.5ex>[r] &\C_0(x)}$, $\xymatrix{\D_1(x) \ar@<.7ex>[r] \ar@<-.7ex>[r] &\D_0(x)}$, respectively.  These are precisely $\C_1(x)$ and $\D_1(x)$, respectively.
\end{rmk}

For our work with reduced $n$-crossed complexes, a truncated version of the above lemma will be useful; the corollary is below and follows from showing injectivity in degree $n$ which we leave to the reader.

\begin{pro} \label{crsisomorphism}
A morphism $f:[\C:\p]\rightarrow [\D:\d]$ of $n$-crossed complexes is a trivial fibration if and only if it satisfies the properties of Lemma \ref{BGfibration} and the last map is an isomorphism when $k=n$.
\end{pro}

%For two reduced crossed complexes $[\H:\p]$ and $[\G:\d]$, the \it{product crossed complex} $[\H\times \G:\p\times \d]$ defined by $(\H\times \G)_n=\H_n\times \G_n$ with differentials being the usual product is a crossed complex.  Moreover, $\H\times \G$  with degree wise projections $\pi_1:\H\times \G \rightarrow \H$ and $\pi_2:\H\times \G \rightarrow \G$ is universal up to isomorphism in $\xc$.

%\begin{Def}
%A reduced crossed complex $\G$ is \it{acyclic} if all homotopy groups are trivial.
%\end{Def}

\begin{pro}
Let $[\H:\p]$ and $[\G:\d]$ be reduced crossed complexes where $\G$ is acyclic. The induced morphism $\pi_1:\H\times \G \rightarrow \H$ is a trivial fibration where $\H\times \G$ is the product crossed complex.
\end{pro}
\begin{proof}
The differentials of the crossed complex $\H\times \G$ are the product maps.  Thus, $\pi_1(\H\times \G)=\coker(\p\times \d) \cong \coker \p \times \coker \d =\pi_1(\H)\times \pi_1(\G)$ and $\pi_k(\H\times \G)=H(\H\times \G) \cong H(\H)_k \times H(\G)_k=\pi_k(\H) \times \pi_k(\G)$ for $k\geq 2$.  Since $\G$ is acyclic, the right side of these products are trivial.  The fact that $\pi_1:\H\times \G\rightarrow \H$ is a fibration is immediate.  
\end{proof}

The following theorem is an extension of the Dold-Kan correspondence and will be useful later.

\begin{AshleyExt} \cite{Ashley2} \label{extension}
There is a fully faithful functor $\Xc \rightarrow \sgrpd $ which restricts to a fully faithful functor $\xc \rightarrow\sgrp$ where $\sgrpd$ and $\sgrp$ are the categories of simplicial groupoids and simplicial groups, respectively.
\end{AshleyExt}

\section{Derived Mapping Groupoid of Crossed Complexes}   

The set of morphisms from a crossed complex $\H$ to a crossed complex $\G$ in the homotopy category of crossed complexes with respect to the Brown-Golasi\'{n}ski model structure (see \cite{BrownGolasinski1}) is the set
\begin{equation}\label{homotopyVSquotient}
[\H,\G]_{\Xc} ~\cong ~ \Xc(Q\H,R\G)/\simeq ~\cong ~\Xc(Q\H,\G)/\simeq
\end{equation}
where $R$ and $Q$ are the fibrant and cofibrant replacement functors, respectively, and the last isomorphism follows from the fact that all crossed complexes are fibrant.  We would like a groupoid $\cG$ such that $[\H,\G]_{\Xc}\cong \pi_0(\cG)$ which we now define beginning with the objects.

\begin{Def}
Let $\H$ and $\G$ be crossed complexes, and let $\Q$ be a cofibrant replacement of $\H$.  A \it{derived morphism} from $\H$ to $\G$ is a morphism in $\Xc(\Q,\G)$.
\end{Def}  

Derived morphisms are conveniently represented as fractions of the form:
\begin{equation}\label{known}
\xymatrix@=1em{ &\Q \ar[dr]^f \ar@{->>}[dl]_p^\simeq \\ \H & &\G}
\end{equation}
where $p$ is a cofibrant replacement of $\H$.  The set of derived morphisms from $\H$ to $\G$ is called the \it{derived mapping set} and denoted by $Rhom(\H,\G)$.  Explicitly, $Rhom(\H,\G)=\Xc(\Q,\G)_0$ where $\Q$ is a cofibrant replacement of $\H$.  

For derived morphisms $f$ and $g$, a morphism from $f$ to $g$ is a homotopy $h$ from $f$ to $g$.  Diagramatically, a morphism is of the form
\begin{equation}
\xymatrix@=1em{&&\Q\ar[dd]|{1_\Q}  \ar@{->>}[dll]^\simeq_p \ar[drr]^f \ar@{}[drr]^(.5){}="b"  &&\\ \H &&&&\G  \\ &&\Q \ar@{->>}[ull]^p_\simeq \ar[urr]_g  \ar@{}[uur]^(.15){}="a" \ar@{<=}"a";"b"_h }
\end{equation}
where $\Theta$ is the identity map.  We will say that the right-side of the diagram commutes up to the homotopy $h$.  

For computational ease, we would prefer the morphisms to be defined algebraically.  Fortunately, Andrew Tonks showed that a homotopy is equivalent to a $1$-fold left homotopy in \ref{eqofhomotopytheories} which is reminiscent of a chain homotopy between morphisms of chain complexes.  Due to the equivalence, we will interchange freely between the two. 

\begin{Def}
Let $\H$ and $\G$ be crossed complexes.  The \it{derived mapping groupoid} from $\H$ to $\G$ is the groupoid $\XC(\Q,\G)_1$ of the internal crossed complex $\XC(\Q,\G)$ where $\Q$ is a cofibrant replacement of $\H$.   The derived mapping groupoid is denoted by $\rhom(\H,\G)$.
\end{Def} 

The following corollary follows directly from Theorem \ref{eqofhomotopytheories}.

\begin{cor}
Let $H$ and $\G$ be crossed complexes.  Then there exists a bijection
$$[\H,\G]_\Xc\cong \pi_0(\rhom(\H,\G)).$$
\end{cor}

Now we begin our analysis of derived morphisms between reduced crossed complexes.  From a derived morphism, we would like to algebraically construct a crossed complex $\E$ that defines a new fraction as below.

\begin{equation}\label{dia:E}
\xymatrix@=1.5em{ &\Q \ar@{->>}[ddl]|\simeq_p \ar[ddr]^f \ar[d]  \\ &\E  \ar@{->>}[dl]^{p^\ast}|\simeq \ar[dr]_{f^\ast} \\ \H &&\G}
\end{equation}

The crossed complex $\E$ will be weakly equivalent to $\Q$ by the 2-out-of-3 property of model categories and should at least be easier to compute than $\Q$ which will give a more manageable model of the derived morphisms of crossed complexes.  An example of an object $\E$ in \ref{dia:E} is the product $\H\times \G$ with the product map and the projections.  Unfortunately, $\pi_1$ is not necessarily a trivial fibration except when $\G$ is acyclic.  However, the universal property of products guarantees that $\E$ would need to factor $p\times f$.  The $n$-pushout construction below will help us construct such an object.

%%%%%%%%%%%%%%%%%%%%%%% n-pushout %%%%%%%%%%%%%%%%%%%%%%%%%

Recall from \cite{Noohi2} that for a group $L$ and the solid diagram 
\begin{equation} \label{FiberedCoproduct}
\xymatrix{ H\ar[r]^p \ar[d]_d & G\ar@{-->}[d]^{q_N\circ i_2} \\ K \ar@{-->}[r]_{q_N\circ i_1} &K\times^H G }
\end{equation}
of $L$-modules, the fibered coproduct of $K$ and $G$ under $H$ is the $L$-module $K\times^H G = (K\times G)/N$ where $N=\{(d(h)\i,p(h)) | h\in H\}$.  The induced morphisms are the usual inclusions into the product $i_1,i_2$ followed by the quotient map $q_N:K\times G\rightarrow K\times^H G$.

We leave the proof of the next theorem to the reader and note that it does not rely on the groups being abelian.

\begin{pro} \label{coproductkernelcokernel}
For the fibered coproduct in \ref{FiberedCoproduct}, the induced morphism $\ker(d) \rightarrow \ker(q_N\circ i_2$) is a surjection and the induced morphism $\coker(d)\rightarrow \coker(q_N\circ i_2)$ is an isomorphism.
\end{pro}

Consider a morphism of reduced crossed complexes $f:[\H:\p]\rightarrow [\G:\d]$.  For each degree $k\geq 2$, $f_k$ and $\p_k$ form a diagram of $\H_1$-modules as in \ref{FiberedCoproduct}.  The quotient $\H_1\ltimes^{\H_2} \G_2$ was shown to be well-defined in \cite{Noohi2}; we will show that the quotient $\H_{k-1}\ltimes^{\H_k} \G_k$ is defined for $k>2$.

\begin{lem} \label{commutativeimage}
Let $[\G:\d]$ be a reduced crossed complex.  Then $\d_3(\G_3)$ is in the center of $\G_2$.
\end{lem}
\begin{proof}
Let $x\in \G_2$ and $z\in \G_3$.  By property CM2 of crossed modules and the fact that $\d_3$ preserves the action of $\G_1$ on $\G_3$ of crossed complexes, $x\i \p_3(z)\i x = (\p_3(z)\i)^{\p_2(x)}  = \p_3(z^{\p_2(x)})\i = \p_3(z)\i$ where the last equality follows from the crossed complex axiom that $\d(\G_2)$ acts trivially on $\G_3$.
\end{proof}

\begin{pro} \label{pushoutWelldefined}
Let $f:[\H:\p]\rightarrow [\G:\d]$ be a morphism of reduced crossed complexes.  Then the quotient $\H_{k-1}\times^{\H_k} \G_k$ is well-defined for $k\geq 3$; in particular, $\N_2$ is in the center of $\H_2\times \G_3$.
\end{pro}
\begin{proof}
Since $\H_k$ is a $\H_1$-submodule, and $\p_k$ and $f_k$ are homomorphisms invariant under the action of $\H_1$ for $k\geq 2$, $\N_{k-1}$ is a $\H_1$-submodule.
For $k\geq 4$, the groups of the homomorphisms are abelian so the quotient $\H_{k-1}\times^{\H_k} \G_k=(\H_{k-1}\times \G_k)/\N_{k-1}$ is well-defined.  

We now show that $\N_2$ is normal in $\H_2\times \G_3$.   Let $z\in \H_2$ and $(x,y)\in \H_2\times \G_3$.  Then 
\begin{align*}
(x,y)\i (\p_3(z)\i,f_3(z))(x,y) 
&= (x\i,y\i) (\p_3(z)\i,f_3(z))(x,y) \\
&= (x\i \p_3(z)\i x,y\i f_3(z) y) \\
&= (\p_3(z)\i x\i x,f_3(z)y\i y) \\
&= (\p_3(z)\i,f_3(z))
\end{align*}
where the last equality holds since $\p_3(\H_3)$ is in the center of $\H_2$ and $\G_3$ is abelian.
Thus, $N_2$ is normal in $\H_2\times \G_3$ and is in fact in the center of $\H_2\times \G_3$.  Hence, the quotient $\H_{k-1}\times^{\H_k} \G_k$ is well-defined for $k\geq 3$.
\end{proof}

\begin{pro} \label{projectionlift}
Let $f:[\H:\p]\rightarrow [\G:\d]$ be a morphism of reduced crossed complexes.  Then there is a lift of the morphism $\p_{k-1}\circ \pi_1:\H_{k-1}\times \G_k \rightarrow \H_{k-2}$ to $\H_{k-1}\times^{\H_k} \G_k$ denoted by $\ol{\p_{k-1}\circ \pi_1}$ for $k\geq 3$.
\end{pro}
\begin{proof}
For $k\geq 3$, 
%we have the solid diagram 
%\begin{equation*}
%\xymatrix{ \H_{k-1}\times \G_k \ar[dr]_{q_N} \ar[r]^(.6){\pi_1} &\H_{k-1} \ar[r]^{\p_{k-1}} &\H_{k-2}  \\  &\H_{k-1}\times^{\H_k} \G_k \ar@{-->}[ur]_(.6){\ol{\p_{k-1}\circ \pi_1}} }
%\end{equation*}
since $\p_{k-1}(\pi_1((\p_k(x)\i,f_k(x))))=\p_{k-1}(\p_k(x)\i)  =\p_{k-1}(\p_k(x\i)) =1 $ for all $x\in \H_k$, the morphism $\p_{k-1}\circ \pi_1$ lifts to $\ol{\p_{k-1}\circ \pi_1}$.
\end{proof}

\begin{Def} \label{npushout}
Let $f:[\H:\p]\rightarrow [\G:\d]$ be a morphism of reduced $n$-crossed complexes.  Then the \it{$n$-pushout below $f$} is the chain of group homomorphisms
\begin{equation}
\xymatrix{\G_n \ar[r]^(.3){q_n\circ i_2} &\H_{n-1}\times^{\H_n}\G_n \ar[r]^(.6){\ol{\p_{k-1}\circ \pi_1}} &\H_{n-2} \ar[r]^{\p_{n-2}} &\cdots \ar[r]^{\p_3} &\H_2 \ar[r]^{\p_2} &\H_1 }.
\end{equation}
\end{Def} 

\begin{pro}  \label{npushoutcrs}
Let $f:[\H:\p]\rightarrow [\G:\d]$ be a morphism of reduced $n$-crossed complexes.  Then the $n$-pushout below $f$ is a reduced $n$-crossed complex denoted by $[\H^{f_n}:\p^{f_n}]$.
\end{pro}

The propositions above and below are cumbersome, but easily verified.  

\begin{pro} \label{npushoutfactor}
Let $f:[\H:\p]\rightarrow [\G:\d]$ be a morphism of reduced $n$-crossed complexes.  Then $[\H^{f_n}:\p^{f_n}]$ factors $f$ by
\begin{equation}
\xymatrix{ \H \ar@/^1pc/[rr]^f \ar[r]_{\iota} &\H^{f_n} \ar[r]_{\rho} &\G}.
\end{equation}
wher $\iota$ and $\rho$ are morphisms of crossed complexes.
\end{pro}

\begin{pro}  \label{cotruncationsurjection}
Let $f:[\H:\p]\rightarrow [\G:\d]$ be a morphism of reduced $n$-crossed complexes.  The induced morphism $\pi_k(\H)\rightarrow \pi_k(\H^{f_n})$ is a surjection for $k=n$ and an isomorphism for $k<n$.
\end{pro}
\begin{proof}
The morphism $\H\rightarrow \H^{f_n}$ has the form of the solid diagram below.
\begin{equation}\label{n-1cotruncationdiagram}
\xymatrix@R=1em{ \H_n \ar[rrr]^{f_n} \ar[d]_{\p_n} &&&\G_n \ar[d]^{q_N\circ i_2}  \ar[d]    \\  \H_{n-1} \ar@{-->}[dr]^{q_{\p_n}} \ar[dd]_{\p_{n-1}} \ar[rrr]^{q_N\circ i_1}  & &&\H_{n-1}\times^{\H_{n}} \G_n  \ar[dd]^{\ol{\p_{n-1}\circ \pi_1}} \ar@{-->}[dl]_q \\  &\H_{n-1}/\p_n \ar@{-->}[dl]^{\ol{\p}_{n-1}} \ar@{-->}[r]^(.3){c} & (\H_{n-1}\times^{\H_{n}} \G_n)/q_N\circ i_2  \ar@{-->}[dr]^l  \\  \H_{\leq n-2}  \ar@{=}[rrr] &&& \H_{\leq n-2}    }
\end{equation}

Since the top square is just the pushout, $\pi_n(\H)\rightarrow \pi_n(\H^{f_n})$ is surjective by Proposition \ref{coproductkernelcokernel}.  Also, the induced map $c$ of cokernels is an isomorphism which satisfies $c\circ q_{\p_n} = q \circ (q_N\circ i_1)$. 

Since $\H$ and $\H^{f_n}$ are chain complexes, there are lifts $\ol{\p}_{n-1}$ and $l$ which make the complete diagram \ref{n-1cotruncationdiagram} commute. Since $c$ is an isomorphism, $c$ induces inverse maps on the kernels of $\ol{\p}_{n-1}$ and $l$.  Thus, the induced morphism $\ker \ol{\p}_{n-1} \rightarrow \ker l$ is an isomorphism.  Since $\ker \ol{\p}_{n-1} \cong \ker \p_{n-1}/\Im\p_{n+1} =\pi_{n-1}(\H)$ and $\ker l\cong \ker \ol{\p_{n-1}\circ \pi_1}/\Im q_N\circ i_2 =\pi_{n-1}(\H^{f_n})$ the induced morphism $\pi_{n-1}(\H) \rightarrow \pi_{n-1}(\H^{f_n})$ is an isomorphism.

For $k=n-2$, the kernels are in fact equal.  We will show that the images of $\p_{n-1}$ and $\ol{\p_{n-1}\circ \pi_1}$ are also equal.  For $x\in \Im \p_{n-1}$, there is a $z\in \H_{n-1}$ such that $\p_{n-1}(z) =x$.  Then $\ol{\p_{n-1}\circ \pi_1}([z,1])=\p_{n-1}(z)=x$. Thus, $x\in \Im \ol{\p_{n-1}\circ \pi_1}$.  

Conversely, for $x\in \Im \ol{\p_{n-1}\circ \pi_1}$, there exists an element $[u,v]\in \H_{n-1}\times^{\H_{n}}\G_n$ such that $x=\ol{\p_{n-1}\circ \pi_1}([u,v]) =l(q([u,v])) =\ol{\p}_{n-1} ( c\i ( q([u,v]))) $.  Since $ c\i ( q([u,v])) \in \H_{n-1}/\p_n$ and $q_{\p_n}$ is surjective, there exists an element $z\in \H_{n-1}$ such that $q_{\p_n}(z)= c\i ( q([u,v]))$.  
Then $\ol{\p}_{n-1} ( c\i ( q([u,v])))  = \ol{\p}_{n-1} ( q_{\p_n}(z))  = \p_{n-1}(z)$.  Thus, $x\in \Im \p_{n-1}$ and the images are equal.  

Hence, we have the equality $\pi_{n-2}(\H)= \pi_{n-2}(\H^{f_n})$.  For $k<n-2$, clearly we have the equality $\pi_k(\H) = \pi_k(\H^{f_n})$. by definition of $\H^{f_n}$.  Hence, we have the desired result.
\end{proof}

\begin{pro}\label{n-1cotruncation}
Let $f:[\H:\p]\rightarrow [\G:\d]$ be a morphism of reduced $n$-crossed complexes.  The $(n-1)$-cotruncation of the morphism $\iota:\H\rightarrow \H^{f_n}$ is the morphism of reduced $n-1$-crossed complexes 
\begin{equation*}
\xymatrix@=1em{ \H_{n-1}/\p_n \ar[r]^(.4){c} \ar[d]_{\ol{\p_{n-1}}} &\H_{n-1}\times^{\H_n}\G_n/\p^f_n \ar[d]^{\ol{\p^f_{n-1}}} \\  \H_{n-2} \ar@{=}[r] &\H_{n-2}}
\end{equation*}
where $c$ is an isomorphism.
\end{pro}
\begin{proof}
The result follows from the commutativity of diagram \ref{n-1cotruncationdiagram} and Proposition \ref{coproductkernelcokernel}.
\end{proof}

The next section will apply the $n$-pushout of reduced $n$-crossed complexes to a specific morphism which will lead to a model of derived morphisms.

\section{N-Butterflies of Reduced N-Crossed Complexes} %\texorpdfstring{$n$}{n}

We now construct an algebraic object from a derived morphism between reduced $n$-crossed complexes by applying the $n$-pushout construction to the induced diagonal morphism.  Unfolding the $n$-pushout results in a diagram which we call a $n$-butterfly named after the paper \cite{Noohi2}.  In the next section, we will show that certain $n$-butterflies model the morphisms between reduced $n$-crossed complexes of $Ho(\Xc)$.

For any morphisms $p:[\Q:\xi] \rightarrow [\H:\p]$ and $f:[\Q:\xi] \rightarrow [\G:\p]$ in $n\xc$, we will denote the induced diagonal morphism by $\nabla^f:\Q\rightarrow \H\times \G$.  Recall that the $n$-pushout is then denoted by $\Q^{\nabla^f_n}$.

\begin{pro}\label{npushoutweakeq}
Let $p:[\Q:\xi]\rightarrow [\H:\p]$ and $f:[\Q:\xi]\rightarrow [\G:\d]$ be morphisms in $n\xc$ where $p$ is a trivial fibration.  Then the canonical morphism  $\iota: \Q\rightarrow \Q^{\nabla^f_n}$ is a weak equivalence.
\end{pro}

\begin{proof}
By proposition \ref{cotruncationsurjection}, the induced morphism of homotopy groups $\pi_n(\Q)\rightarrow \pi_n(\Q^{\nabla^f_n})$ is an isomorphism for $1\leq k <n$ and a surjection for $k=n$.  We only need to show that the induced morphism above is injective.  

The morphism $\Q \rightarrow \Q^{\nabla^f_n}$ in degrees $n-1$ and $n$ has the form

\begin{equation} 
\xymatrix@=1em{ \ker \xi_n \ar[d] \ar[r]  &\ker i_2 \ar[d]  \\  \Q_n \ar[r]^{\nabla_n} \ar[d]_{\xi_n}  &\H_n \times \G_n \ar[d]^{i_2}  \\  \Q_{n-1} \ar[r]^(.3){i_1} \ar[d] & \Q_{n-1}\times^{\Q_n} (\H_n \times \G_n) \ar[d] \\ \Q_{\leq n-2} \ar@{=}[r] &\Q_{\leq n-2} }
\end{equation}

The surjection of the kernels above is given explicitly by Proposition \ref{coproductkernelcokernel}.  Since $p:\Q\rightarrow \H$ is an acyclic fibration, $\Q_n \cong \ker \xi_{n-1} \times_{\ker \p_{n-1}} \H_n$ by proposition \ref{crsisomorphism}.  Moreover, $\xi_n=\pi_1$ and $p_n=\pi_2$.  Let $(a,b)\in \ker \xi_{n-1} \times_{\ker \p_{n-1}} \H_n$ such that $(a,b)\in \ker \xi_n$ and $\nabla_n( (a,b)) = (1,1)$.
Since $\xi_n=\pi_1$ and $p_n=\pi_2$, $1=\xi_n((a,b)) = a$ and $(1,1) = \nabla_n( (a,b)) = (p_n((a,b)), f_n((a,b))) = (b, f_n((a,b)))$.  Thus, $a=1$ and $b=1$; in other words, $(a,b)=(1,1)$.

We have shown that the induced morphism $\pi_n(\Q) =\ker \xi_n \rightarrow \ker i_2 =\pi_n(\Q^{\nabla^f_n})$ is injective, completing the proof that the homotopy groups in degree $1\leq k\leq n$ are isomorphic.
\end{proof}

\begin{cor}\label{npushoutfraction}
Let $p:[\Q:\xi]\rightarrow [\H:\p]$ and $f:[\Q:\xi]\rightarrow [\G:\d]$ be morphisms in $n\xc$ where $p$ is a trivial fibration.  Then the composition $\xymatrix{\Q^{\nabla^f_n} \ar[r]^\rho &\H\times \G \ar[r]^{\pi_1}  &\H}$ is a trivial fibration.
\end{cor}
\begin{proof}
Since the trivial fibration $p$ factors as $p=\pi_1\circ \rho \circ \iota$, the morphism $\pi_1\circ \rho$ is a weak equivalence by Proposition \ref{npushoutweakeq} and the two out of three property of model categories.  The fibration follows from the definition of fibration of reduced crossed complexes and the set theoretic fact about factorizations of set surjections.
\end{proof}

For the next proposition consider the diagram below induced by the canonical morphism $\rho:\Q^{\nabla^f_n}\rightarrow \H\times \G$.
\begin{equation*}
\xymatrix{ \G_n \ar[r]^(.4){i_2} &\H_n\times \G_n  \ar[d]_(.4){\xi^{\nabla^f_n}_n } &  \\  &\Q^{\nabla^f_n}_{n-1} \ar[d]_{\xi^{\nabla^f_n}_{n-1}} \ar[r]^{\pi_1\circ \rho_{n-1}}  &\H_{n-1} \ar[d]^{\p_{n-1}} \\  &\ker \xi_{n-2} \ar[r]_{\pi_1\circ \rho_{n-2}} &\ker \p_{n-2}}
\end{equation*}

\begin{pro}
Let $p:[\Q:\xi]\rightarrow [\H:\p]$ and $f:[\Q:\xi]\rightarrow [\G:\d]$ be morphisms in $n\xc$ where $p$ is a trivial fibration.  Then the induced chain complexes 
\begin{equation*}
\xymatrix{ 1\rightarrow \G_n \ar[r]^i &\Q^{\nabla^f_n}_{n-1} \ar[r]^(.3)u & \ker \xi_{n-2}\times_{\ker \p_{n-2}} \H_{n-1} \ar[r] &1 }
\end{equation*}
\begin{equation*}
\xymatrix{  \Q^{\nabla^f_n}_k \ar[r]^(.3){u_k} & \ker \xi_{k-1}\times_{\ker \p_{k-1}} \H_k \ar[r] &1  }
\end{equation*}
are exact where $i= \xi^{\nabla^f_n}\circ i_2$ and $k\leq n-2$.
\end{pro}
\begin{proof}
We first prove that $i$ is injective.
Suppose $a\in \G_n$ such that $i(a)=[1,(1,1)]$.
Since $i(a)=[1,(1,a)]$ by definition of $i$, $(1,(1,a))=(\xi_n(x)\i,(p_n(x),f_n(x)))$ for some $x\in \Q_n$.  Since $p:\Q\rightarrow \H$ is a weak equivalence, $p_n$ induces the isomorphism $\ker \xi_n\cong \ker \p_n$.  Since $x\in \ker \xi_n$ and $p_n(x)=1$, $x=1$.  Since $f_n(x)=a$ and $f_n(x)=1$, $a=1$.  Hence, the sequence considered is exact at $\G_n$.

Suppose $[a,(b,c)]\in \ker u$.  Then 
\begin{align*}
(1,1)
=u([a,(b,c)]) &=(\xi^{\nabla^f_n}_{n-1}([a,(b,c)]),\pi_1\circ \rho_{n-1}([a,(b,c)])) \\
&=(\xi_{n-1}\circ \pi_1([a,(b,c)]),\pi_1\circ \rho_{n-1}([a,(b,c)])) \\
&= (\xi_{n-1}(a),\pi_1((p_{n-1}(a)\p_n(b),f_{n-1}(a)\d_n(x)))).
\end{align*}
Since $p_{n-1}(a)=\p_n(b\i)$, $\pi_1\circ \rho_{n-1}([a,(1,1)])
=\pi_1\circ \rho_{n-1}(\iota_{n-1}(a))  
=\pi_1\circ (\rho_{n-1}\circ \iota_{n-1})(a)  
=\pi_1((p_{n-1}(a),f_{n-1}(a))) 
=p_{n-1}(a)  
=\p_n(b\i)$.

Since $a\in \ker \xi_{n-1}$, $[a,(1,1)]\in \ker \xi^{\nabla^f_n}_{n-1}$ and $[a,(1,1)]$ represents an element in the homology group $H(\Q^{\nabla^f_n})_{n-1}$ which we will denote by $[[a,(1,1)]]$.
By definition of the induced morphism on homology $(\pi_1\circ \rho_{n-1})^\ast:H(\Q^{\nabla^f_n})_{n-1}\rightarrow H(\H)_{n-1}$, $(\pi_1\circ \rho_{n-1})^\ast([[a,(1,1)]])=[\pi_1\circ \rho_{n-1}([a,(1,1)])]=[\p_n(b\i)]=[1]$ in $H(\H)_{n-1}$.

Since $\pi_1\circ \rho$ is a trivial fibration by Proposition \ref{npushoutfraction}, the induced morphism on homology groups is in fact an isomorphism.  Thus, $[[a,(1,1)]]=[ [1,(1,1)]]$ in $H(\Q^{\nabla^f_n})_{n-1}$ so there exists $(x,y)\in \H_n\times \G_n$ such that $[a,(1,1)]=\xi^{\nabla^f_n}_n((x,y))=[1,(x,y)]$.  Then $[a,(b,c)]
= [a,(1,1)][1,(b,c)] 
= [1,(x,y)][1,(b,c)] 
= [1,(xb,yc)]$.

Since $\p_n(x)=\pi_1\circ \rho_{n-1}([1,(x,y)]) =\pi_1\circ \rho_{n-1}([a,(1,1)]) =p_{n-1}(a)$, $\p_n(xb)=\p_n(x)\p_n(b)=p_{n-1}(a)\p_n(b)=1$.  Since $xb\in \ker \p_n$ and $\ker \xi_n=H(\Q)_n\cong H(\H)_n = \ker \p_n$, there exists a $z\in \ker \xi_n$ such that $p_n(z)=xb$.  
Thus, 
\begin{align*}
(1,(xb,yc))&=(\xi_n(z)\i,(p_n(z),yc f_n(z\i)f_n(z))) \\
&=(1,(1,ycf(z\i)))(\xi_n(z)\i,(p_n(z),f_n(z))).
\end{align*}

By definition of $\Q^{\nabla^f_n}_{n-1}$, $[a,(b,c)]=[1,(xb,yc)]=[1,(1,ycf(z\i))]$ where $ycf(z\i)\in \G_n$.
Hence, $i(ycf(z\i))=[a,(b,c)]$ in $\Q^{\nabla^f_n}_{n-1}$ and the homology is trivial at $\Q^{\nabla^f_n}_{n-1}$.

Since $\pi_1\circ \rho:\Q^{\nabla^f_n}\rightarrow \H$ is a trivial fibration, the induced morphisms:
\begin{equation*}
\xymatrix{  \Q^{\nabla^f_n}_k \ar[r]^(.3){u_k} & \ker \xi_{k-1}\times_{\ker \p_{k-1}} \H_k }
\end{equation*}
are surjections when $k\leq n-1$ by Proposition \ref{crsisomorphism}.
\newline Hence, the complete complex of the proposition is exact.
\end{proof}

The above analysis in this section has characterized the $n$-pushout below $\nabla^f:\Q\rightarrow \H\times \G$.  We know unfold this $n$-pushout and abstract its characteristics to obtain $n$-butterflies.

\begin{Def}
Let $[\H:\p],[\G:\d]$ be reduced $n$-crossed complexes.  A \it{$n$-butterfly} $([\E,\eta],p,f,\alpha, \beta)$ from $\H$ to $\G$ is a commutative diagram in Fig. \ref{fig:butterfly} and satisfying the axioms below.
%\begin{equation*}
%\xymatrix@=0.5em{ \H_n \ar[dr]^\alpha \ar[dd]_{\p_n} & &\G_n \ar[dd]^{\d_n} \ar[dl]_\beta \\ &\E_{n-1} \ar[dr]^{f_n} \ar[dl]_{p_n} \ar[dd]^{\eta_{n-1}} \\ \H_{n-1} \ar[dd]_{\p_{n-1}} && \G_{n-1} \ar[dd]^{\d_{n-1}} \\ &\E_{\leq n-2} \ar[dr]^f \ar[dl]_p \\ \H_{\leq n-2} && \G_{\leq n-2} }
%\end{equation*}
\begin{enumerate}
\item[$\B1$] the maps $[\E:\eta]\rightarrow [\H_{\leq n-1}:\p ]$ and $[\E:\eta]\rightarrow [\G_{\leq n-1}:\p]$ are morphisms of reduced $(n-1)$-crossed complexes;
\item[$\B2$]
the induced sequences $\xymatrix{ 1\ar[r]& \G_n \ar[r]^\beta &\E_{n-1} \ar[r]^(.25){u_{n-1}} & \ker \eta_{n-2}\times_{\ker \p_{n-2}} \H_{n-1} \ar[r] &1  }$ and $\xymatrix{  \E_k \ar[r]^(.25){u_k} & \ker \eta_{k-1}\times_{\ker \p_{k-1}} \H_k \ar[r] &1  }$ are exact for $k\leq n-2$;
\item[$\B3$] the compositions $\eta_{n-1}\circ (\alpha\times\beta)$ and $f_n\circ \alpha$ are complexes;
\item[$\B4$] $\alpha,\beta$ satisfy the compatibility conditions $\alpha(x^{p_1(a)})=\alpha(x)^a$ and $\beta(y^{f_1(a)})=\beta(y)^a$ where $x\in \H_n$ and $y\in \G_n$.
\end{enumerate}
\end{Def}

\begin{figure}
\centering
\begin{minipage}{.4\textwidth}
  \centering
\begin{equation*}
\xymatrix@=0.5em{ \H_n \ar[dr]^\alpha \ar[dd]_{\p_n} & &\G_n \ar[dd]^{\d_n} \ar[dl]_\beta \\ &\E_{n-1} \ar[dr]^{f_n} \ar[dl]_{p_n} \ar[dd]^{\eta_{n-1}} \\ \H_{n-1} \ar[dd]_{\p_{n-1}} && \G_{n-1} \ar[dd]^{\d_{n-1}} \\ &\E_{\leq n-2} \ar[dr]^f \ar[dl]_p \\ \H_{\leq n-2} && \G_{\leq n-2} }
\end{equation*}
  \caption{An $n$-Butterfly}%{A figure}
  \label{fig:butterfly}
\end{minipage}
\begin{minipage}{.5\textwidth}
  \centering
\begin{equation*}
\xymatrix@=0.4em{ \H_n \ar[dr] \ar[dd] & &\G_n \ar[dd] \ar[dl] \\ &\Q_{n-1} \times^{\nabla^f_n} (\H_n\times \G_n) \ar[dr] \ar[dl] \ar[dd] \\ \H_{n-1} \ar[dd] && \G_{n-1} \ar[dd] \\ &\Q_{\leq n-2} \ar[dr] \ar[dl] \\ \H_{\leq n-2} &&\G_{\leq n-2} }
\end{equation*}  
\caption{$n$-Butterfly of Pushout}%{Another figure}
  \label{fig:pushoutbutterfly}
\end{minipage}
\end{figure}
The fact that $\eta_{n-1}\circ (\alpha\times\beta)$ is the trivial homomorphism and the homomorphisms $\alpha,\beta$ satisfy the compatibility conditions guarantees the induced complex 
\begin{equation}\label{fig:foldedbutterfly}
\xymatrix{\H_n\times \G_n \ar[r]^{\alpha \times \beta} &\E_{n-1} \ar[r]^{\eta_{n-1}} &\E_{n-2} \ar[r] &\cdots \ar[r]&\E_1  }
\end{equation}
is a reduced $n$-crossed complex which we will denote by $[\E^*:\eta^*]$ with natural maps $p^*:\E^*\rightarrow \H$ and $f^*:\E^*\rightarrow \G$.

\begin{ex}
Let $f$ be an object of $n\xc(\Q,\G)$.  Unfolding the induced morphism $\Q^{\nabla^f_n}\rightarrow \H \times \G$ of reduced $n$-crossed complexes defines a $n$-butterfly.  Explicitly, we have the diagram in \ref{fig:pushoutbutterfly}.
%\begin{equation}\label{inducedbutterfly}
%\xymatrix@=0.4em{ \H_n \ar[dr] \ar[dd] & &\G_n \ar[dd] \ar[dl] \\ &\Q_{n-1} \times^{\nabla^f_n} (\H_n\times \G_n) \ar[dr] \ar[dl] \ar[dd] \\ \H_{n-1} \ar[dd] && \G_{n-1} \ar[dd] \\ &\Q_{\leq n-2} \ar[dr] \ar[dl] \\ \H_{\leq n-2} &&\G_{\leq n-2} }
%\end{equation}
\end{ex}

For the remainder of the paper, we will denote the induced butterfly by $\B^f$ and the middle reduced $(n-1)$-crossed complex by $[\ol{\Q}^{\nabla^f_n}:\ol{\xi}^{\nabla^f_n}]$.  Since $\Q^{\nabla^f_n}$ and $\ol{\Q}^{\nabla^f_n}$ are cumbersome notations, we will use $\Q^f$ and $\ol{\Q}^f$, respectively, as long as they are clear in context i.e. do not confuse with the pushout below $f$.  In summary, 
\begin{equation*}
\B^f=([\ol{\Q}^f:\ol{\xi}^f],\pi_1\circ \rho^f_,\pi_2\circ \rho^f,\xi^f_n\circ i_1, \xi^f_n\circ i_2).
\end{equation*}

\begin{pro}  \label{butterflytrivialfibration}
Let $[\H,\p],[\G,\d]$ be reduced $n$-crossed complexes and $([\E,\eta],p,f,\alpha, \beta)$ be a $n$-butterfly from $\H$ to $\G$.  Then the induced morphism $[\E^*,\xi^*]\rightarrow [\H,\p]$ of reduced $n$-crossed complexes is a trivial fibration. 
\end{pro}
\begin{proof}
By axiom $\B2$ of $n$-butterflies and Proposition \ref{crsisomorphism}, we only need to prove that the induced morphism $\xymatrix{ \H_n\times \G_n \ar[r]^(.35)k &\ker \eta_{n-1} \times_{\ker \p_{n-1}} \H_n  }$ is an isomorphism.
Suppose $(a,b)\in \H_n\times \G_n$ such that $k((a,b))=(1,1)$.  Then $(1,1)=k((a,b))=(\alpha(a)\beta(b),a)$.  Since $\beta(b)=\alpha(1)\beta(b)=\alpha(a)\beta(b)=1$, $b\in ker \beta$.  
By the axiom $\B2$ of $n$-butterflies, $\beta$ is injective.  Thus, $b=1$ and $(a,b)=(1,1)$.  Hence, $k$ is injective.

Suppose $(a,b)\in \ker \eta_{n-1} \times_{\ker \p_{n-1}} \H_n$.  Then $\eta_{n-1}(a)=1$ and $p_{n-1}(a)=\p_n(b)$.  Since $p_{n-1}(a)=\p_n(b)=\p_n(\pi_1((b,1)))=p_{n-1}(\alpha\times \beta(b,1))=p_{n-1}(\alpha(b)\beta(1))=p_{n-1}(\alpha(b))$, $p_{n-1}(a\alpha(b)\i)=p_{n-1}(a)p_{n-1}(\alpha(a))\i=1$.  Then the induced morphism
\begin{equation*}
\xymatrix{ \E_{n-1} \ar[r]^(.25){u_{n-1}} &\ker \eta_{n-2}\times_{\ker \p_{n-2}} \H_{n-1} }
\end{equation*}
satisfies $u_{n-1}(a\alpha(b)\i) 
= (\eta_{n-1}(a\alpha(b)\i),p_{n-1}(a\alpha(b)\i))  
=(\eta_{n-1}(a)\eta_{n-1}(\alpha(b\i)),1)  
=(\eta_{n-1}(\alpha\times \beta ((b\i,1))),1)$.

By axiom $\B3$ of $n$-butterflies, $u_{n-1}(a\alpha(b)\i)=(1,1).  $By axiom $\B2$ of $n$-butterflies, there exists an element $z\in \G_n$ such that $\beta(z)=a\alpha(b)\i$.  Thus, $k((b,z))
=(\alpha\times \beta((b,z)),b) 
=(\alpha(b)\beta(z),b)  
=(\alpha(b)a\alpha(b)\i,b)  
=(a\alpha(a)\alpha(b)\i,b) 
=(a,b)$.
In the case when $n=3$, the commutativity applied above is valid by Proposition \ref{commutativeimage}
Hence, $k$ is surjective and indeed an isomorphism.
\end{proof}

\begin{cor}
Let $[\H,\p],[\G,\d]$ be reduced $n$-crossed complexes and $([\E,\eta],p,f,\alpha, \beta)$ be a $n$-butterfly from $\H$ to $\G$.  If $[\Q:\xi]$ is a cofibrant replacement of $[\H:\p]$, then there exists a weak equivalence $l:\Q\rightarrow \E^*$ such that the following diagram commutes.
\begin{equation} \label{firstbutterflylift}
\xymatrix{ \ast \ar[r] \ar[d] &\E^* \ar@{->>}[d]^{p^*}|\simeq   \\  \Q \ar@{-->}[ur]^l \ar@{->>}[r]^\simeq_p &\H }
\end{equation}
\end{cor}
\begin{proof}
By proposition \ref{butterflytrivialfibration}, $p^*:\E^*\rightarrow \H$ is a trivial fibration. 
Since $\Q$ is a cofibrant object, we have a lift $l$ in \ref{firstbutterflylift}.  By the two-out-of-three axiom of model categories $l$ is a weak equivalence.
\end{proof}

\begin{rmk} \label{aboutQbutterflies}
Our goal is to have a completely algebraic ``butterfly'' which models derived morphisms as in the case for  crossed modules; however, the $n$-pushout construction reduces the cofibrant object by only one degree.  In particular, Proposition \ref{n-1cotruncation} shows that the induced morphism $cotr_{n-1}(\Q)\rightarrow cotr_{n-1}(\Q^f)$ is an isomorphism in degree $n-1$ and the identity in all other degrees.  We now adapt this characteristic to the $n$-butterflies and we hope later work will completely eliminate the topological part of $n$-butterflies.
\end{rmk}

\begin{Def}
Let $\H,\G$ be reduced $n$-crossed complexes and $\Q$ be a cofibrant replacement of $\H$.  A \textit{$n$-butterfly over $\Q$} is a $n$-butterfly $([\E,\eta],p,f,\alpha, \beta)$ such that there exists a lift $l$ in \ref{firstbutterflylift} which induces a morphism $cotr_{n-1}(\Q)\rightarrow cotr_{n-1}(\E^*)$ that is an isomorphism in degree $n-1$ and the identity in all other degrees.  
\end{Def}

\begin{Def}
A morphism of $n$-butterflies over $\Q$ from $\H$ to $\G$
\begin{equation*}
\xymatrix{  (\Theta,\phi):([\E,\eta],p,f,\alpha, \beta) \ar[r] &([\E\me,\eta\me],p\me,f\me,\alpha\me, \beta\me )  }
\end{equation*}
is a morphism of reduced $(n-1)$-crossed complexes $\Theta:\E \rightarrow \E'$ such that $\Theta_{n-1}$ is a group isomorphism, $\Theta_{k}=1_{\E_k}$ for $k\leq n-2$, and Fig. \ref{fig:butterflymorphism}
%\begin{equation*}
%\xymatrix@=0.5em{  &&&\G_n \ar[ddl] \ar[ddd] \ar[dll] \\ \H_n \ar[r] \ar[drr] \ar[ddd] &\E_{n-1}\ar[dr]^{\Theta_{n-1}} \ar@{.>}[ddd] \ar@{.>}@/^.6pc/[ddrr] \ar@{.>}[dddl] \\ &&\E\me_{n-1} \ar[ddd]\ar[dr] \ar[ddll] \\   && &\G_{n-1} \ar[ddd]  \\ \H_{n-1} \ar[ddd] &\E_{\leq n-2} \ar@{.>}@/^.6pc/[ddrr] \ar@{=}[dr] \ar@{..>}[lddd] \\ &&\E\me_{\leq n-2}  \ar[dr] \ar[ddll] \\ &&&\G_{\leq n-2}  \\ \H_{\leq n-2} }
%\end{equation*}
commutes up to a $1$-fold left homotopy (Fig. \ref{fig:butterflyhomotopy})
%\begin{equation*}
%\xymatrix{  \E \ar[rr]^f \ar[dr]_\Theta &\ar@{=>}[d]^\phi &\G  \\  &\E\me \ar[ur]_{f\me} }
%\end{equation*}
of morphisms between reduced $(n-1)$-crossed complexes.
\end{Def}

\begin{figure}
\centering
\begin{minipage}{.5\textwidth}
  \centering
\begin{equation*}
\xymatrix@=0.5em{  &&&\G_n \ar[ddl] \ar[ddd] \ar[dll] \\ \H_n \ar[r] \ar[drr] \ar[ddd] &\E_{n-1}\ar[dr]^{\Theta_{n-1}} \ar@{.>}[ddd] \ar@{.>}@/^.6pc/[ddrr] \ar@{.>}[dddl] \\ &&\E\me_{n-1} \ar[ddd]\ar[dr] \ar[ddll] \\   && &\G_{n-1} \ar[ddd]  \\ \H_{n-1} \ar[ddd] &\E_{\leq n-2} \ar@{.>}@/^.6pc/[ddrr] \ar@{=}[dr] \ar@{..>}[lddd] \\ &&\E\me_{\leq n-2}  \ar[dr] \ar[ddll] \\ &&&\G_{\leq n-2}  \\ \H_{\leq n-2} }
\end{equation*}
  \caption{Morphism of $n$-Butterflies}%{A figure}
  \label{fig:butterflymorphism}
\end{minipage}
\begin{minipage}{.3\textwidth}
  \centering
\begin{equation*}
\xymatrix{  \E \ar[rr]^f \ar[dr]_\Theta &\ar@{=>}[d]^\phi &\G  \\  &\E\me \ar[ur]_{f\me} }
\end{equation*}
  \caption{Homotopy}%{Another figure}
  \label{fig:butterflyhomotopy}
\end{minipage}
\end{figure}

Let $\B=(\E,p,f,\alpha,\beta)$ be an $n$-butterfly over $\Q$.  The identity morphism $1_\B:\B\rightarrow \B$ is the pair $(1_\E,\phi^{id})$ where $1_\E$ is the identity morphism of $\E$ and $\phi^{id}$ is the constant pointed $1$-fold left homotopy defined by $\phi^{id}(x)=1$.  

Now we consider composition of morphisms of $n$-butterflies.  For $n=2$ (the butterfly case), a morphism is an isomorphism $\Theta$ that commutes between butterflies and composition is simply composition of isomorphisms.  For $n>2$, the composition of two morphisms $(\Theta^2,\phi^2), (\Theta^1,\phi^1)$ is defined to be $(\Theta^2,\phi^2)\circ (\Theta^1,\phi^1)=(\Theta^2\circ \Theta^1,\phi^2 \ast \phi^1)$ where $\phi^2*\phi^1(x)=\phi^2(x)\phi^1(x)$.  In fact, choosing the morphism to be $\Theta^2\circ \Theta^1$ completely determines the homotopy by equation \ref{pointeddetermine}.

\begin{pro}
Let $[\H,\p],[\G,\d]$ be reduced $n$-crossed complexes and $\Q$ be a cofibrant replacement of $\H$.  The $n$-butterflies from $\H$ to $\G$ over $\Q$ and their morphisms form a groupoid denoted by $n\B_\Q(\H,\G)$.  
\end{pro}
\begin{proof}
The identity and associative properties of composition of $n$-butterfly morphisms is clear from the definition of $\phi^{id}$ and $\phi^2*\phi^1$.  Let $(\Theta,\phi):([\E^u:\eta^u],p^u,u)\rightarrow ([\E^v:\eta^v],p^v,v)$ be a morphism of $n$-butterflies.  Then $p=p'\circ \Theta$.  Define $\Theta\i:\E^v\rightarrow \E^u$ in degree $n-1$ as the inverse of the isomorphism $\Theta_{n-1}$ and the identity in degrees less than $n-1$.  Clearly, $\Theta\i \circ \Theta=1_{\E^u}$ and $\Theta \circ \Theta\i=1_{\E^v}$.  Since $p'(\Theta(x))=p(x)$, $p(\Theta\i(x))=p'(\Theta(\Theta\i(x)))=p'(x)$.  

Defining a $1$-fold left homotopy $\phi\i$ by $\phi\i(x)=\phi(x)\i$ gives a morphism a morphism of $n$-butterflies $(\Theta\i,\phi\i):([\E^v:\eta^v],p^v,v)\rightarrow ([\E^u:\eta^u],p^u,u)$.  Moreover, $(\Theta\i,\phi\i)\circ (\Theta,\phi)=(\Theta\i\circ \Theta,\phi\i*\phi)=(1_{\E^u},\phi^{id})$ and $(\Theta,\phi)\circ (\Theta\i,\phi\i)=(\Theta\circ \Theta\i,\phi*\phi\i)=(1_{\E^v},\phi^{id})$ which satisfies the properties of the inverse of $(\Theta,\phi)$. 
%Since homotopies are invertible, we can define a homotopy $\lambda$ by the diagram
%\begin{equation*}
%\xymatrix{  \E' \ar[rr]^{f'} \ar[dr]_{1_{\E'}} \ar[dd]_{\Theta\i} &\ar@{=>}[d]^{1_g} &\G  \\  &\E' \ar@{=>}[d]^{\phi\i} \ar[ur]_{f'}  \\  \E \ar[ur]^\Theta  \ar@/_3pc/[uurr]_f  &&}
%\end{equation*}
%The morphism $(\Theta\i,\lambda)$ is the inverse of $(\Theta,\phi)$.
\end{proof}

\begin{cor}\label{nbutterflyweakeq}
Let $[\H,\p],[\G,\d]$ be reduced $n$-crossed complexes and 
\begin{equation*}
\xymatrix{  (\Theta,h):([\E,\eta],p,f,\alpha, \beta) \ar[r] &([\E\me,\eta\me],p\me,f\me,\alpha\me, \beta\me )  }
\end{equation*}
be a morphism of $n$-butterflies.
Then the induced morphism $\E^*\rightarrow (\E\me)^*$ of reduced $n$-crossed complexes is a weak equivalence.
\end{cor}
\begin{proof}
By Proposition \ref{butterflytrivialfibration} and two out of three property of weak equivalences.
\end{proof}

\section{Modeling Morphisms of Homotopy N-Types} \label{equivalencesection}

We now show that the groupoid of $n$-butterflies over $\Q$ from $\H$ to $\G$ is equivalent to the groupoid of derived morphisms from $\H$ to $\G$.  Note that we only need to use the condition of being over $\Q$ when proving essentially surjective.  

\begin{thm} \label{equivalencetheorem}
Let $[\H:\p]$ and $[\G:\d]$ be reduced $n$-crossed complexes.  There is an equivalence of groupoids $n\rhom(\H,\G)\simeq n\B_\Q(\H,\G)$. 
\end{thm}

Before we prove there is an equivalence, we show that as a result of the theorem above $n$-butterflies over $\Q$ may be viewed as a partially algebraic model of certain pointed connected homotopy $n$-types. The Dold-Kan Extension \ref{extension} and the fact that simplicial groups model pointed topological spaces define a full embedding $Ho(\xc)\rightarrow Ho(\Top^*)$ where $\Top^*$ is the category of pointed connected homotopy types.  Consequently, for pointed connected homotopy $n$-types $X,Y$ corresponding to reduced $n$-crossed complexes $\H$ and $\G$ via the embedding, $[X,Y]_{\Top^*} \cong [\H,\G]_\xc$.  We now state our main theorem

\begin{thm}\label{mainthm}
Let $[\H:\p]$ and $[\G:\d]$ be reduced $n$-crossed complexes.  Then we have the bijection $[\H,\G]_{\Xc} \cong \pi_0(n\B_\Q(\H,\G))$.
\end{thm}
\begin{proof}
Since homotopy is equivelent to $1$-fold left homotopies $[\H,\G]_{\Xc} \cong \pi_0(\rhom(\H,\G))$.  By the Theorem \ref{equivalencetheorem}, $\pi_0(\rhom(\H,\G))\cong \pi_0(n\B_\Q(\H,\G))$.
\end{proof}

Hence, there is a subclass of pointed connected homotopy $n$-types which are modeled by the category of reduced $n$-crossed complexes where the morphisms are the connected components of the groupoid of $n$-butterflies over $\Q$.  However, we restricted our $n$-butterflies to those with a special lifting property to have essential surjectivity.  We hope future work will produce a purely algebraic model.

We know begin proving Theorem \ref{equivalencetheorem} by first realizing there exists a functor from the groupoid of derived morphisms to $n$-butterflies.  The proof will follow a lemma.

\begin{pro}\label{omegafunctor}
Let $[\H:\p]$ and $[\G:\d]$ be reduced $n$-crossed complexes.  There is a functor $\Omega: n\Xc(\Q,\G)\rightarrow n\B(\H,\G)$ where $p:[\Q:\xi]\rightarrow [\H:\p]$ is a cofibrant replacement of $\H$.
\end{pro}

\begin{lem} \label{inducedbutterflymorphism}
Let $f,g$ be objects in $n\Xc(\Q,\H)$ and $h:f\Rightarrow g$ be a homotopy.  Then there exists an induced morphism of $n$-butterflies $(\Theta,h^\ast):\B^f\rightarrow \B^g$.
\end{lem}
\begin{proof}
Since $h:f\Rightarrow g$ is a homotopy of morphisms of reduced $n$-crossed complexes, there is a pointed $1$-fold left homotopy $(g,\phi_k:\Q_k\rightarrow \G_{k+1})$ satisfying \ref{pointeddetermine} for $k< n$ and $f_n(x)=g_n(x)\phi_{n-1}(\xi_n(x))$.  We will first define the morphism $\Theta:\ol{\Q}^f\rightarrow \ol{\Q}^g$ of reduced $(n-1)$-crossed complexes.  Define $\Theta_{n-1}:\Q^f_{n-1}\rightarrow \Q^g_{n-1}$ by $\Theta_{n-1}([a,(b,c)]_f)=[a,(b,\phi_{n-1}(a)c)]_g$.  By the Preservation Law of $\phi$, $\Theta_{n-1}$ is a group homomorphism for $n\geq 2$.

Since both of the crossed complexes $\Q^f$ and $\Q^g$ are weakly equivalent to $\H$ and they are equal in degrees $k\leq n-2$, we have the isomorphisms in the following diagram.  
\begin{equation} \label{isomorphismtheta}
\xymatrix@=1.5em{ \ker \xi^f_n \ar[r] \ar[d]|\cong &\H_n\times \G_n \ar[r]^{\xi^f_n} \ar@{=}[d] &\Q^f_{n-1} \ar[r] \ar[d]^{\Theta_{n-1}} &\coker \xi^f_n \ar[r] \ar[d]|\cong &\1 \ar@{=}[d] \\  \ker \xi^g_n \ar[r] &\H_n\times \G_n \ar[r]^{\xi^g_n} &\Q^g_{n-1} \ar[r] &\coker \xi^g_n \ar[r] &\1}
\end{equation}
Clearly, the left and right squares of \ref{isomorphismtheta} commute.
Since $\phi_{n-1}(1)=1$ for $n\geq 2$ by the preservation laws of $\phi$, 
$\Theta_{n-1}(\xi^f_n((a,b)))
=\Theta_{n-1}([1,(a,b)]_f) 
=[1,(a,\phi_{n-1}(1)b)]_g 
=[1,(a,b)]_g  
=\xi^g_n((a,b))$
for $(a,b)\in \H_n\times \G_n$.
Thus, the middle-left square of \ref{isomorphismtheta} commutes.  
The commutativity of the middle-right square follows from the fact that $(1,\phi_{n-1}(a))\in \H_n\times \G_n$ for any $a\in \Q_{n-1}$.  By the five lemma, $\Theta_{n-1}$ is a group isomorphism.

Letting $\Theta_k=1_{\Q_k}$ for $1\leq k\leq n-2$ completes the morphism $\Theta:\ol{\Q}^f \rightarrow \ol{\Q}^g$ of reduced $(n-1)$-crossed complexes.  Moreover, the left side of the following commutes. 
\begin{equation} \label{ThetaAboveView}
\xymatrix@=1em{&&\ol{\Q}^f\ar[dd]|{\Theta}  \ar[dll]_{\pi_1\circ \rho^f} \ar[drr]^{\pi_2\circ \rho^f} \ar@{}[drr]^(.5){}="b"  &&\\ \H_{<n} &&&&\G_{<n}  \\ &&\ol{\Q}^g \ar[ull]^{\pi_1\circ \rho^g} \ar[urr]_{\pi_2\circ \rho^g}  \ar@{}[uur]^(.15){}="a" \ar@{<=}"a";"b"_{h^*} }
\end{equation}

We now claim that defining $\phi_k$ for $1\leq k\leq n-2$ makes the right-side of Diagram \ref{ThetaAboveView} commute up to homotopy.  Since $\ol{\Q}^f,\ol{\Q}^g$ are reduced $(n-1)$-crossed complexes and $\Q^f_k=\Q_k=\Q^g_k$ for $1\leq k \leq n-2$, we only need to check that $\phi$ satisfies the pointed $1$-fold homotopy condition in degree $n-1$.  Explicitly, $\pi_2\circ \rho^f_{n-1}([a,(b,c)])=$
\begin{align}
(\pi_2\circ \rho^g_1)(\Theta_1([a,(b,c)])) \t{  for  } n=2 \\
(\pi_2\circ \rho^g_{n-1})( \Theta_{n-1}([a,(b,c)])\phi_{n-2}(\xi^f_{n-1}([a,(b,c)]))  \t{ for } n>2 \label{thetahomotopy}
\end{align}
The butterfly case (n=2) has already been proven.  For $n>2$, the left-side of the equality is $\pi_2\circ \rho^f_{n-1}([a,(b,c)]) = \pi_2((p_{n-1}(a),f_{n-1}(a))(\p_n(b),\d_n(c))) =f_{n-1}(a)\d_n(c)$.  We now calculate the right-side of equation \ref{thetahomotopy} for $n=3$ and $n=4$ (see App. \ref{nequals3functorproof} and \ref{nequals4functorproof}, respectively).  Hence, we have the desired equality for $n\geq 2$ and furthermore, a induced morphism $h^\ast=(\Theta,\phi):\B^f\rightarrow \B^g$ of $n$-butterflies.
\end{proof}

We now finish the proof that there exists a functor from derived morphisms to $n$-butterflies.

\begin{proof}[Proof of Theorem \ref{omegafunctor}]
Define $\xymatrix{ \Omega:n\Xc(\Q,\G)_0\rightarrow n\B(\H,\G)_0  }$ by assigning to a morphism $f:\Q\rightarrow \G$ in $n\Xc(\Q,\G)$ the $n$-butterfly $\Omega(f)$ obtained by unfolding the morphism $\Q^{\nabla^f_n}\rightarrow \H\times \G$ of reduced $n$-crossed complexes as in Fig. \ref{fig:foldedbutterfly}. %\ref{inducedbutterfly}.  

Suppose $f,g\in ob(n\Xc(\Q,\G))$.  Define $\xymatrix{ \Omega:n\Xc(\Q,\G)(f,g) \rightarrow n\B(\H,\G)(\Omega(f),\Omega(g))  }$ by assigning to a homotopy $h:f\Rightarrow g$ the induced morphism of $n$-butterflies $h^\ast=(\Theta,\phi):\Omega(f)\rightarrow \Omega(g)$ from Lemma \ref{inducedbutterflymorphism}.  We only need to show that $\Omega(1_f)=1_{\Omega(f)}$ and $\Omega(h_2\circ h_1) = \Omega(h_2)\circ \Omega(h_1)$.

The identity homotopy of $f$ in $n\Xc(\Q,\G)$, $1_f:f\Rightarrow f$, is defined by the constant pointed $1$-fold left homotopy $\phi_{id}(x)=1$.  By Lemma \ref{inducedbutterflymorphism}, the induced morphism $\Omega(1_f)$ is defined as $(\Theta,\phi_{id})$ where $\Theta_{n-1}([a,(b,c)]) = [a,(b,\phi_{id}(a)c)] = [a,(b,1c)]=[a,(b,c)]$ and $\Theta_k=1_{\Q_k}$ for $1\leq k \leq n-2$.  Thus, $\Omega(1_f)=(1_\Q,\phi_{id})=1_{\Omega(f)}$.

Consider two morphisms $h_1:u\Rightarrow v$ and $h_2:v\Rightarrow w$ in $n\Xc(\Q,\G)$ defined by the pointed $1$-fold left homotopies $\phi^1,\phi^2$, respectively.  By definition of $\Omega$, $\Omega(h_1)=(\Theta^1,\phi^1)$, $\Omega(h_2)=(\Theta^2,\phi^2)$ where $\Theta^1$ and $\Theta^2$ are defined by $\Theta^1([a,(b,c)])=[a,(b,\phi^1(a)c)]$ and $\Theta^2([a,(b,c)])=[a,(b,\phi^2(a)c)]$, respectively.  

By definition of composition in $n\B(\H,\G)$, $\Omega(h_2)\circ \Omega(h_1)=(\Theta^2,\phi^2)\circ (\Theta^1,\phi^1) = (\Theta^2\circ \Theta^1, \phi^2*\phi^1)$ where $\Theta^2_{n-1}\circ \Theta^1_{n-1}([a,(b,c)])=\Theta^2_{n-1}([a,(b,\phi_1(a)c)]) = [a,(b,\phi_2(a)\phi_1(a)c)]$. Since $h_2\circ h_1$ is defined by the $1$-fold left homotopy $\phi^2*\phi^1(x)=\phi^2(x)\phi^1(x)$, $\Omega(h_2\circ h_1)=(\Theta^{2,1},\phi^2*\phi^1)$ where $\Theta^{2,1}_{n-1}$ is defined by $\Theta^{2,1}_{n-1}([a,(b,c)])=[a,(b,\phi^2*\phi^1(a)c)]= [a,(b,\phi_2(a)\phi_1(a)c)]=\Theta^2_{n-1}\circ \Theta^1_{n-1}([a,(b,c)])$. 

Hence, $\Omega(h_2\circ h_1) =  (\Theta^{2,1},\phi^2*\phi^1)=(\Theta^2\circ \Theta^1,\phi_2*\phi_1)= (\Theta^2,\phi_2)\circ (\Theta^1,\phi_1) = \Omega(h_2)\circ \Omega(h_1)$.
\end{proof}

\begin{pro}
Let $[\H:\p]$ and $[\G:\d]$ be reduced $n$-crossed complexes and $\Q$ be a cofibrant replacement of $\H$.  The functor $\Omega:n\Xc(\Q,\G)\rightarrow n\B(\H,\G)$ is faithfull. 
\end{pro}
\begin{proof}
Since $n=2$ is the butterfly case, we are only concerned with $n>2$.  Let $f,g:\Q\rightarrow \G$ be morphisms of reduced $n$-crossed complexes.  We must show that the set map $n\Xc(\Q,\G)(f,g)\rightarrow n\B(\H,\G)(\Omega(f),\Omega(g))$ is injective.
Suppose $h_1,h_2:f\Rightarrow g$ are two homotopies determined by the $1$-fold left homotopies $(g,\phi_k), (g,\lambda_k)$ and the induced morphisms of $n$-butterflies $\Omega(h_1)=(\Theta^\phi,\phi^*_k)$ and $\Omega(h_2)=(\Theta^\lambda,\lambda^*_k)$ from $\Omega(f)$ to $\Omega(g)$ are equal.
Then $\Theta^\phi=\Theta^\lambda$ and $\phi^*_k=\lambda^*_k$ for all $k<n-1$.  Since $\phi^*_k=\phi_k$ and $\lambda^*_k=\lambda_k$ for $k<n-1$ by definition, we only need to show that $\phi_{n-1}=\lambda_{n-1}$ to prove that $h_1=h_2$.

\sloppy Let $x\in \Q_{n-1}$.  Since $\Theta^\phi=\Theta^\lambda$, $[x,(1,\phi_{n-1}(x))]_g=\Theta^\phi_{n-1}([x,(1,1)]_f)=\Theta^\lambda_{n-1}([x,(1,1)]_f)=[x,(1,\lambda_{n-1}(x))]_g$ in $\Q^g_{n-1}$.  By definition of equality in $\Q^g_{n-1}$, there exists $z\in \Q_n$ such that for $n>2$, $(\xi_n(z)\i,\nabla^g_n(z))=(x,(1,\phi_{n-1}(x)))\i(x,(1,\lambda_{n-1}(x))) =(1,(1,\phi(x)\i_{n-1}\lambda_{n-1}(x)))$.  Then $z\in \ker \xi_n$, $p_n(z)=1$ and $g_n(z)=\phi(x)\i_{n-1}\lambda_{n-1}(x)$ for any $n\geq 2$.  The first two observations imply $z=1$ since $p:\Q\rightarrow \H$ is a trivial fibration; in particular, $p_n$ induces an isomorphism $\ker \xi_n \cong \ker \p_n$.  Thus, $\phi(x)\i_{n-1}\lambda_{n-1}(x)=g_n(1)=1$; equivalently, $\phi_{n-1}(x)=\lambda_{n-1}(x)$.  Hence, $\phi_{n-1}=\lambda_{n-1}$ and we have proven injectivity.
\end{proof}

\begin{pro}
Let $[\H:\p]$ and $[\G:\d]$ be reduced $n$-crossed complexes and $\Q$ be a cofibrant replacement of $\H$.  The functor $\Omega:n\Xc(\Q,\G)\rightarrow n\B(\H,\G)$ is full. 
\end{pro}

\begin{proof}
Let $f,g:\Q\rightarrow \G$ be morphisms of reduced $n$-crossed complexes.  We must show that the set map $n\Xc(\Q,\G)(f,g)\rightarrow n\B(\H,\G)(\Omega(f),\Omega(g))$ is surjective.  Suppose that there is a morphism of butterflies $(\Theta,h):\Omega(f)\rightarrow \Omega(g)$.  The homotopy $h$ induces a pointed $1$-fold left homotopy of reduced $(n-1)$-crossed complexes given by
\begin{equation*}
(g\circ \Theta, \phi_k:\Q^f_k\rightarrow \G_{k+1}) \text{~for~} 1\leq k\leq n-1
\end{equation*}
Furthermore, $f$ is determined by 
\begin{align}
f_1(x)&=g_1(x)\d_2(\phi_1(x));  \\
f_k(x)&=g_k(x)\d_{k+1}(\phi_k(x))\phi_{k-1}(\xi_k(x)) \text{~for~}~1\leq k <n-1;  \\
(\pi_2\circ \rho^f)_{n-1}([a,(b,c)])&=((\pi_2\circ \rho^g)_{n-1}\circ \Theta_{n-1})([a,(b,c)])\phi_{n-2}(\xi^f_{n-1}([a,(b,c)])). \label{pirhothetahomotopy}
\end{align}

We will define a homotopy $h^\ast$ from $f$ to $g$, or equivalently, a pointed $1$-fold left homotopy $(g,\phi^\ast_k:\Q_k\rightarrow \G_{k+1}) \text{~for~} 1\leq k\leq n-1$ such that 
\begin{align} 
f_1(x)&=g_1(x)\d_2(\phi^\ast_1(x)); \label{injhomotopy1} \\
f_k(x)&=g_k(x)\d_{k+1}(\phi^\ast_k(x))\phi^\ast_{k-1}(\xi_k(x)) \text{~for~}~1\leq k\leq n-1;  \label{injhomotopy2} \\
f_{n}(x)&=g_n(x)\phi^\ast_{n-1}(\xi_n(x)). \label{injhomotopy3}
\end{align}
Since $\Theta_k:\Q^{\nabla^f_n}\rightarrow \Q^{\nabla^g_n}_k$ is the identity map for $k<n-1$, we let $\phi^\ast_k=\phi_k$ for $k<n-1$.   By definition of $\phi$, \ref{injhomotopy1} and \ref{injhomotopy2} for $k<n-1$ are true.  We define $\phi^\ast_{n-1}:\Q_{n-1}\rightarrow \G_n$ by definining a map $\varphi_{n-1}:\Q_{n-1}\rightarrow \H_n\times \G_n/\nabla^g_n(\ker \xi_n)$ which lifts to $\H_n\times \G_n$ and then projecting to $\G_n$.  Informally, for $q\in \Q_{n-1}$ we take $\varphi_{n-1}(q)$ to be the difference between $\Theta_{n-1}([q,(1,1)]_f)$ and $[q,(1,1)]_g$ in $\Q^g$.  

Let $[a,(b,c)]_g=\Theta([q,(1,1)]_f)$ in $\Q_{n-1}\times^{\nabla^g_n}(\H_n\times \G_n)$.  By proposition \ref{coproductkernelcokernel} and the universal property of cokernels on $\Theta$, $[a]=[q]$ in $\coker \xi_n$.  Thus, $a=q\xi_n(x)$ for some $x\in \Q_n$ and 
\begin{align*}
\Theta_{n-1}([q,(1,1)]_f)&=[a,(b,c)]_g  \\
&=[q\xi_n(x),(b,c)]_g  \\
&=[q\xi_n(x),(b,c)\nabla^g_n(x)\nabla^g_n(x)\i]_g  \\
%&=[q,(b,c)\nabla^g_n(x)]_g[\xi_n(x),\nabla^g_n(x)\i]_g \\
&=[q,(b,c)\nabla^g_n(x)]_g[\xi_n(x\i)\i,\nabla^g_n(x\i)]_g \\
%&=[q,(bp_n(x),cg_n(x))]_g [1]_g \\
%&=[q,(bp_n(x),cg_n(x))]_g  \\
&=[q,(1,1)]_g[1,(bp_n(x),cg_n(x))]_g
\end{align*}
Further, since $\pi_1\circ \rho^f_{n-1}=(\pi_1\circ \rho^g_{n-1})\circ \Theta_{n-1}$ by definition of the morphism $(\Theta,h)$, we have $\pi_1(\rho^f_{n-1}([q,(1,1)]_f)=\pi_1(\nabla^f_{n-1}(q)(\p_n(1),\d_n(1))) =p_{n-1}(q)$ and the following. 
\begin{align*}
%\pi_1(\rho^f_{n-1}([q,(1,1)]_f)&=\pi_1(\nabla^f_{n-1}(q)(\p_n(1),\d_n(1))) \\
%&=\pi_1((p_{n-1}(q),f_{n-1}(q))(1,1))  \\
%&=p_{n-1}(q)  \\
\pi_1(\rho^g_{n-1}(\Theta_{n-1}([q,(1,1)]_f)))&=\pi_1(\rho^g_{n-1}([q,(1,1)]_g[1,(bp_n(x),cg_n(x))]_g))  \\
%&=\pi_1(\rho^g_{n-1}([q,(bp_n(x),cg_n(x))]_g))  \\
&=\pi_1(\nabla^g_{n-1}(q)(\p_n(bp_n(x)),\d_n(cg_n(x))))  \\
%&=\pi_1((p_{n-1}(q),g_{n-1}(q))(\p_n(bp_n(x)),\d_n(cg_n(x))))  \\
%&=\pi_1((p_{n-1}(q)\p_n(bp_n(x)),g_{n-1}(q)\d_n(cg_n(x)))  \\
&=p_{n-1}(q)\p_n(bp_n(x))
\end{align*}
Thus, $\p_n(bp_n(x))=1$ and $bp_n(x)\in \ker \p_n$.  Since $\ker \p_n\cong \ker \xi_n$, there exists $z\in \ker \xi_n$ such that $p_n(z)=bp_n(x)$.  Let $M=\nabla^g_n(\ker \xi_n)$.  Then $[bp_n(x),cg_n(x)]_M = [bp_n(x),cg_n(x)]_M[p_n(z),g_n(z)]\i_M =[1,cg_n(xz\i)]_M$ in $\H_n\times \G_n/M$.  

Define $\varphi_{n-1}:\Q_{n-1}\rightarrow \H_n\times \G_n/M$ by $\varphi_{n-1}(q)=[1,cg_n(xz\i)]_M$.  Suppose $(a',(b',c'))$ is another representative of $\Theta([q,(1,1)]_f)$ i.e. $\Theta([q,(1,1)]_f)=[a',(b',c')]_g$.
Similarly to $a$, there exists a $x'\in \Q_n$ such that $a'=q\xi_n(x')$ and $\Theta([q,(1,1)]_f)=[q,(1,1)]_g[1,(b' p_n(x'),c' g_n(x'))]_g$.  There exists $z' \in \ker \xi_n$ such that $[b' p_n(x'),c' g_n(x')]_M = [1,c' g_n(x'(z')\i)]_M$.
Thus, $\varphi_{n-1}(q)=[1,c' g_n(x'(z')\i)]_M$.  Since $[q,(bp_n(x),cg_n(x))]_g=\Theta([q,(1,1)]_f)=[q,(b' p_n(x'),c' g_n(x'))]_g$ in $\Q^g_{n-1}$, there exists a $w\in \Q_n$ such that 
\begin{align*}
(\xi_n(w)\i,\nabla^g_n(w))&=(q,(bp_n(x),cg_n(x)))(q,(b' p_n(x'),c' g_n(x')))\i  \\
&=(1,(bp_n(x),cg_n(x))(b' p_n(x'),c' g_n(x'))\i)
\end{align*}
Then $\nabla^g_n(w)=(bp_n(x),cg_n(x))(b' p_n(x'),c' g_n(x'))\i$ where $w\in \ker \xi_n$.  Thus, 
\begin{equation*}
[1,cg_n(xz\i)]_M=[bp_n(x),cg_n(x)]=[b' p_n(x'),c' g_n(x')]=[1,c'g_n(x'(z')\i)]_M
\end{equation*}
in $\H_n\times \G_n/\nabla^g_n(\ker \xi_n)$.  Hence, $\varphi_{n-1}$ is well-defined.  

{\sloppy Suppose $s,t\in \Q_{n-1}$ and $\Theta_{n-1}([st,(1,1)]_f)=[a,(b,c)]_g$, $\Theta_{n-1}([s,(1,1)]_f) =  [a_s,(b_s,c_s)]_g$, and $\Theta_{n-1}([t,(1,1)]_f)=[a_t,(b_t,c_t)]_g$.  }Since $\Theta_{n-1}$ is an isomorphism, 
\begin{align*}
[st,(1,1)]_g[1,(bp_n(x),cq_n(x))]_g&=\Theta_{n-1}([st,(1,1)]_f) \\
&=\Theta_{n-1}([s,(1,1)]_f)\Theta_{n-1}([t,(1,1)]_f)  \\
%&=[s,(1,1)]_g[1,(b_sp_n(x_s),c_sg_n(x_s))]_g[t,(1,1)]_g[1,(b_tp_n(x_t),c_tg_n(x_t))]_g  \\
&=[st,(1,1)]_g[1,(b_sp_n(x_s)b_tp_n(x_t),c_sg_n(x_s)c_tg_n(x_t))]_g
\end{align*}
where $x, x_s, s_t\in \Q_n$ such that $a=q\xi(x), a_s=q\xi_n(x_s), a_t=q\xi_n(x_t)$.
By cancelling on the left, $[1,(bp_n(x),cg_n(x))]_g=[1,(b_sp_n(x_s)b_tp_n(x_t),c_sg_n(x_s)c_tg_n(x_t))]_g$.  Then there exists $w\in \Q_n$ such that 
\begin{equation*}
(\xi_n(w)\i,\nabla^g_n(w))(1,(bp_n(x),cg_n(x)))=(1,(b_sp_n(x_s)b_tp_n(x_t),c_sg_n(x_s)c_tg_n(x_t))).
\end{equation*}
In particular, $[bp_n(x),cg_n(x)]_M=[b_sp_n(x_s),c_sg_n(x_s)]_M=[b_tp_n(x_t),c_tg_n(x_t)]_M$.
\begin{align*}
\varphi_{n-1}(st)&=[1,cg_n(xz\i)]_M  \\
&=[bp_n(x),cg_n(x)]_M  \\
&=[b_sp_n(x_s),c_sg_n(x_s)]_M[b_tp_n(x_t),c_tg_n(x_t)]_M  \\
&=[1,c_sg_n(x_sz_s\i)]_M[1,c_tg_n(x_tz_t\i)]_M  \\
&=\varphi_{n-1}(s)\varphi_{n-1}(t)
\end{align*}
Hence, $\varphi_{n-1}$ is a homomorphism.

Since $\Q\rightarrow \Q^g$ is a weak equivalence of reduced $n$-crossed complexes, $\nabla^g_n(\ker \xi_n) = \ker \xi^g_n$.  Then the natural morphism $\C\rightarrow \D$ of reduced crossed complexes given by
\begin{equation*}
\xymatrix@=1.5em{ \cdots \ar[r]  &\1 \ar[r] \ar[d] &\ker \xi^g_n \ar[d]\ar[r]^i &\H_n\times \G_n \ar[r] \ar[d]^{q_\nabla} &\1 \ar[r] \ar[d] &\cdots \ar[r]  &\1 \ar[r] \ar[d] &\Q_1 \ar@{=}[d] \\
 \cdots \ar[r] &\1 \ar[r] &\1 \ar[r] &\H_n\times\G_n/\nabla^g_n(\ker \xi_n) \ar[r] &\1 \ar[r] &\cdots \ar[r] &\1 \ar[r] &\Q_1  }
\end{equation*}
where $\H_n\times \G_n$ is in degree $n-1$ is a trivial fibration.  Since $\Q$ is cofibrant, there exists a lift $\varphi^*:\Q\rightarrow \C$ where $\varphi_{n-1} = q_\nabla \circ \varphi^*_{n-1}$.  Thus, we have the following commutative diagram.
\begin{equation}  \label{factorphi}
\xymatrix@=1.5em{ \Q_n \ar[dd]_{\xi_n} \ar[r]^{\nabla^g_n} & \H_n\times \G_n\ar[d]^{q_{\nabla}}  \ar[r]^{\pi_2}  &\G_n \ar[dd]^{\d_n} \\ 
&\H_n\times \G_n/\nabla^g_n(\ker \xi_n)  \ar[d]^{\ol{\xi^g_n}}\\ 
\Q_{n-1} \ar@{..>}[ur]^{\varphi_{n-1}} \ar@{..>}@/^1pc/[uur]^{\varphi^*_{n-1}} \ar[r]_{\iota_{n-1}~~~~~~~} &\Q_{n-1}\times^{\nabla^g_n} \H_n\times \G_n \ar[r]_{~~~~~~~\pi_2\circ \rho^g} &\G_{n-1}  }
\end{equation}

Define $\phi^\ast_{n-1}$ to be the composition $\pi_2 \circ \varphi^*_{n-1}:\Q_{n-1}\rightarrow \H_n\times \G_n \rightarrow \G_n$.  We must show that $(g,\phi^\ast_k)$ satisfies equation \ref{injhomotopy2} for $k=n-1$ and \ref{injhomotopy3}.  For $k=n-1$, see App. \ref{kequalsnminus1Full}. 

For $k=n$, let $x\in \Q_n$.  From the morphism of $n$-butterflies, $\xi^g_n=\Theta_{n-1}\circ \xi^f_n$.  Then $\Theta_{n-1}([1,(p_n(x),f_n(x))]_f) 
=\Theta_{n-1}(\xi^f_n(p_n(x),f_n(x))) 
=\xi^g_n(p_n(x)f_n(x))  
=[1,(p_n(x),f_n(x))]_g$.   
By definition of $\phi^*_{n-1}$,
\begin{align*}
\Theta_{n-1}([\xi_n(x),(1,1)]_f) &=[\xi_n(x),(1,1)]_g [1,(1,\phi^*_{n-1}(\xi_n(x)))]_g  \\
&=[1,(p_n(x),g_n(x))]_g [1,(1,\phi*_{n-1}(\xi_n(x)))]_g  
\end{align*}

Since $[1,(p_n(x),f_n(x))]_f=[\xi_n(x),(1,1)]_f$ in $Q^f_{n-1}$, 
\begin{equation*}
[1,(1,f_n(x))]_g=[1,(1,g_n(x)\phi^*_{n-1}(\xi_n(x)))]_g.
\end{equation*}
Thus, there exists $w\in \Q_n$ such that $w\in \ker \xi_n$ and 
\begin{equation*}
(p_n(w),g_n(w))=(1,f_n(x))\i(1,g_n(x)\phi^*_{n-1}(\xi_n(x))).
\end{equation*}
Then $[1,g_n(x)\i f_n(x)]_M=[1,\phi^*_{n-1}(\xi_n(x))]_M$.  By definition of $\phi^*_{n-1}$,  $f_n(x) = \\ g_n(x)\phi^*_{n-1}(\xi_n(x))$.

Lastly, we confirm that $\Omega(\phi^*)=(\Theta,\phi)$.  By definition of $\phi^*_k$, $\phi^*_k=\phi_k$ for $k<n-1$.  By definition of $\Omega$, $\Omega(\phi^*) = (\Psi,\phi)$ where $\Psi_k=1_\Q$ for $k<n-1$ and $\Psi_{n-1}:\Q_{n-1}\times^f \H_n\times \G_n \rightarrow \Q_{n-1}\times^g \H_n\times \G_n$ is defined by $\Psi_{n-1}([a,(b,c)]_f)=[a,(b,\phi^*_{n-1}(a)c)]_g$.  We only need to show that $\Psi_{n-1}=\Theta_{n-1}$ (see App. \ref{psiequalsthetaFull}).  Hence, $(g,\phi^*_k)$ is the desired pointed $1$-fold left homotopy from $f$ to $g$.

\end{proof}

\begin{pro}
Let $[\H:\p]$ and $[\G:\d]$ be reduced $n$-crossed complexes and $\Q$ be a cofibrant replacement of $\H$.  The functor $\Omega:n\Xc(\Q,\G)\rightarrow n\B_\Q(\H,\G)$ is essentially surjective. 
\end{pro}

\begin{proof}
Suppose $\B=([\E:\eta],p,g,\alpha, \beta)$ is a $n$-butterfly over $\Q$.  Folding the butterfly gives us the reduced $n$-crossed complex of \ref{fig:foldedbutterfly} with morphisms $p^*:\E^*\rightarrow \H$ and $g^*:\E^*\rightarrow \H$.  Moreover, $p^*$ is a trivial fibration by proposition \ref{butterflytrivialfibration}.  Since $\B$ is a $n$-butterfly over $\Q$, there exists a lift $l$ in \ref{firstbutterflylift} which induces a morphism $cotr_{n-1}(\Q)\rightarrow cotr_{n-1}(\E^*)$ that is an isomorphism in degree $n-1$ and the identity in all other degrees.  Denote the composition $g^*\circ l: \Q\rightarrow \G$ by $f$ which is a derived morphism with $p$.

We claim that there exists an isomorphism $(\Theta,\phi):\Omega(f)=\B^f \rightarrow \B$.\footnote{Recall that $\B^f=([\ol{\Q}^f:\ol{\xi}^f],\pi_1\circ \rho^f_{\leq n-1},\pi_2\circ \rho^f_{\leq n-1},\xi^f_n\circ i_1,\xi^f_n\circ i_2)$.}  
By definition of a $n$-butterfly over $\Q$, $l_k$ is the identity map for $k<n-1$; therefore, we let $\Theta_k=1_\Q$ for $k< n-1$.  For $k=n-1$, define $l^*_{n-1}$ by $(a,(b,c)) \mapsto l_{n-1}(a)\alpha\times \beta((b,c))$.  We know that $l^*_{n-1}$ is a homomorphism for the buttefly case ($n=2$) and since crossed complexes are abelian for $k\geq 3$, it is a homomorphism for $n\geq 4$.  For $n=3$,
\begin{align*}
l^*_{n-1}((a,(b,c))(a\me,(b\me,c\me)))%&=l^*_{n-1}((aa\me,(bb\me,cc\me)))  \\
&=l_{n-1}(aa\me)\alpha\times \beta((bb\me,cc\me))  \\
%&=l_{n-1}(a)l_{n-1}(a\me)\alpha\times \beta((b,c))\alpha\times\beta((b\me,c\me))  \\
&=l_{n-1}(a)\alpha\times \beta((b,c))l_{n-1}(a\me)\alpha\times\beta((b\me,c\me)) &\text{[Prop. \ref{commutativeimage}]}  \\
&=l^*_{n-1}((a,(b,c))l^*_{n-1}((a\me,(b\me,c\me)))
\end{align*}
To show that $l^*_{n-1}$ lifts to $\Q^f$, let $z\in \Q_n$.  By definition of $l^*$, $p_n=\pi_1\circ l_n$ and $g^*\circ l_n = \pi_2\circ l_n^*$.  Thus, by the universal property of products, $l*_n=p_n\times (g^*\circ l_n) =p_n\times f_n =\nabla^f_n$.  Thus, $l^*_{n-1}((\xi_n(z)\i,\nabla^f_n(z)))=l_{n-1}(\xi_n(z))\i\alpha\times\beta(\nabla^f_n(z)) =\alpha\times\beta(\nabla^f_n(z))\i\alpha\times\beta(\nabla^f_n(z))=1_{\E_n}$.

Denoting the lift by $\Theta_{n-1}$, $(\Theta,\phi^{id})$ is the desired morphism of $n$-butterflies as long as $\Theta_{n-1}$ is an isomorphism.  Proposition \ref{nbutterflyweakeq} implies the induced morphism $\Q^f\rightarrow \E^*$ is a trivial fibration.  The five-lemma along with the commutative diagram of short exact sequences
\begin{equation*}
\xymatrix@R=1.5em{ \1 \ar[r] \ar[d] &\ker \xi^f_n \ar[d]|\cong \ar[r] & \H_n\times \G_n \ar[r] \ar@{=}[d]&\Q^f_{n-1} \ar[d]|{\Theta_{n-1}}\ar[r] &\coker \xi^f_n \ar[r] \ar[d]|\cong &\1 \ar[d]  \\  \1 \ar[r] &\ker \alpha\times \beta \ar[r] & \H_n\times \G_n \ar[r] &\E_{n-1} \ar[r] &\coker \alpha\times \beta \ar[r] &\1 }
\end{equation*}
implies that $\Theta_{n-1}$ is an isomorphism. 
\end{proof}

\appendix
\section{Equations}

Below are lengthy computations from Section \ref{equivalencesection}.

\begin{align} \label{nequals3functorproof}
\begin{split}
&(\pi_2\circ \rho^g_2)\circ \Theta_2([a,(b,c)])\phi_1(\xi^f_2([a,(b,c)])) \\ 
&= \pi_2\circ \rho^g_2([a,(b,\phi_2(a)c])\phi_1(\xi_2(a))  \\
&=\pi_2((p_2(a),g_2(a))(\p_3(b),\d_3(\phi_2(a)c)))\phi_1(\xi_2(a))  \\
&=\pi_2((p_2(a)\p_3(b),g_2(a)\d_3(\phi_2(a)c)))\phi_1(\xi_2(a))  \\
&=g_2(a)\d_3(\phi_2(a)c)\phi_1(\xi_2(a))  \\
&=g_2(a)\d_3(\phi_2(a))\d_3(c)\phi_1(\xi_2(a))  \\
&=g_2(a)\d_3(\phi_2(a))\phi_1(\xi_2(a))\phi_1(\xi_2(a))\i\d_3(c)\phi_1(\xi_2(a))  \\
&=f_2(a)\phi_1(\xi_2(a))\i\d_3(c)\phi_1(\xi_2(a))  \\
&=f_2(a)\d_3(c)^{\d_2((\phi_1(\xi_2(a)))}  \\
&=f_2(a)\d_3(c^{\d_2((\phi_1(\xi_2(a)))})  \\
&=f_2(a)\d_3(c)
\end{split}
\end{align}

The last equality follows from the fact that $\d_2$ acts trivially on $\G_2$.  For $n\geq 4$:

\begin{align} \label{nequals4functorproof}
\begin{split}
&(\pi_2\circ \rho^g_{n-1})\circ \Theta_{n-1}([a,(b,c)])\phi_{n-2}(\xi^f_{n-1}([a,(b,c)])) \\ 
&= \pi_2\circ \rho^g_{n-1}([a,(b,\phi_{n-1}(a)c])\phi_{n-2}(\xi_{n-1}(a))  \\
&=\pi_2((p_{n-1}(a),g_{n-1}(a))(\p_n(b),\d_n(\phi_{n-1}(a)c)))\phi_{n-2}(\xi_{n-1}(a))  \\
&=\pi_2((p_{n-1}(a)\p_n(b),g_{n-1}(a)\d_n(\phi_{n-1}(a)c)))\phi_{n-2}(\xi_{n-1}(a))  \\
&=g_{n-1}(a)\d_n(\phi_{n-1}(a)c)\phi_{n-2}(\xi_{n-1}(a))  \\
&=g_{n-1}(a)\d_n(\phi_{n-1}(a))\d_n(c)\phi_{n-2}(\xi_{n-1}(a))  \\
&=g_{n-1}(a)\d_n(\phi_{n-1}(a))\phi_{n-2}(\xi_{n-1}(a))\d_n(c)  \\
&=f_{n-1}(a)\d_n(c)
\end{split}
\end{align}

\begin{align} \label{kequalsnminus1Full}
\begin{split}
&g_{n-1}(q)\d_n(\phi^*_{n-1}(q))\phi^*_{n-2}(\xi_{n-1}(q)) \\
%&= g_{n-1}(q)\d_n(\phi^*_{n-1}(q))\phi_{n-2}(\xi_{n-1}(q))  \\
&=\pi_2\circ \rho^g(\iota^g(q))\d_n(\phi^*_{n-1}(q))\phi_{n-2}(\xi_{n-1}(q))  \\
%&=\pi_2\circ \rho^g([q,(1,1)]_g)\d_n(\phi^*_{n-1}(q))\phi_{n-2}(\xi_{n-1}(q))  \\
&=\pi_2\circ \rho^g([q,(1,1)]_g[1,(bp_n(x),cg_n(x))]_g[1,(bp_n(x),cg_n(x))]\i_g)\d_n(\phi^*_{n-1}(q))\phi_{n-2}(\xi_{n-1}(q))  \\ 
&=\pi_2\circ \rho^g( \Theta_{n-1}([q,(1,1)]_f)[1,(bp_n(x),cg_n(x))]\i_g )\d_n(\phi^*_{n-1}(q))\phi_{n-2}(\xi_{n-1}(q))   \\
&=\pi_2\circ \rho^g( \Theta_{n-1}([q,(1,1)]_f)) \pi_2\circ \rho^g( [1,(bp_n(x),cg_n(x))]\i_g ) \d_n(\phi^*_{n-1}(q))\phi_{n-2}(\xi_{n-1}(q))  \\
%&=\pi_2\circ \rho^g( \Theta_{n-1}([q,(1,1)]_f)) \pi_2\circ \rho^g( [1,(bp_n(x),cg_n(x))]_g )\i  \d_n(\phi^*_{n-1}(q))\phi_{n-2}(\xi_{n-1}(q))  \\
&=\pi_2\circ \rho^g( \Theta_{n-1}([q,(1,1)]_f)) \pi_2\circ \rho^g(\xi^g_n((bp_n(x),cg_n(x))) )\i  \d_n(\phi^*_{n-1}(q))\phi_{n-2}(\xi_{n-1}(q))  \\
&=\pi_2\circ \rho^g( \Theta_{n-1}([q,(1,1)]_f)) \pi_2\circ \rho^g(\ol{\xi^g_n}([bp_n(x),cg_n(x)]_M) )\i  \d_n(\phi^*_{n-1}(q))\phi_{n-2}(\xi_{n-1}(q))  \\
&=\pi_2\circ \rho^g( \Theta_{n-1}([q,(1,1)]_f)) \pi_2\circ \rho^g(\ol{\xi^g_n}([1,cg_n(xz\i)]_M) )\i  \d_n(\phi^*_{n-1}(q))\phi_{n-2}(\xi_{n-1}(q))  \\
&=\pi_2\circ \rho^g( \Theta_{n-1}([q,(1,1)]_f)) \pi_2\circ \rho^g(\ol{\xi^g_n}(\varphi_{n-1}(q)) )\i  \d_n(\phi^*_{n-1}(q))\phi_{n-2}(\xi_{n-1}(q))  \\
&=\pi_2\circ \rho^g( \Theta_{n-1}([q,(1,1)]_f)) \pi_2\circ \rho^g(\ol{\xi^g_n}(q_\nabla(\varphi^*_{n-1}(q)]_M)) )\i  \d_n(\phi^*_{n-1}(q))\phi_{n-2}(\xi_{n-1}(q))  \\
&=\pi_2\circ \rho^g( \Theta_{n-1}([q,(1,1)]_f)) \d_n(\pi_2(\varphi^*_{n-1}(q)))\i  \d_n(\phi^*_{n-1}(q))\phi_{n-2}(\xi_{n-1}(q)) \text{~~~~~[\ref{factorphi}]}\\
&=\pi_2\circ \rho^g( \Theta_{n-1}([q,(1,1)]_f)) \d_n(\phi^*_{n-1}(q))\i  \d_n(\phi^*_{n-1}(q))\phi_{n-2}(\xi_{n-1}(q))  \\
%&=\pi_2\circ \rho^g( \Theta_{n-1}([q,(1,1)]_f))\phi_{n-2}(\xi_{n-1}(q))  \\
&=\pi_2\circ \rho^g( \Theta_{n-1}([q,(1,1)]_f))\phi_{n-2}(\xi^f_{n-1}([q,(1,1)]))  \\
&=(\pi_2\circ \rho^f)_{n-1}([q,(1,1)])  \text{~~~~~[\ref{pirhothetahomotopy}]} \\
&=f_{n-1}(q) 
\end{split}
\end{align}

\begin{align} \label{psiequalsthetaFull}
\begin{split}
\Psi_{n-1}([q,(u,v)]_f)&=[q,(u,\phi^*_{n-1}(q)v)]_g  \\
&=[q,(1,1)]_g[1,(1,\phi^*_{n-1}(q)]_g[1,(u,v)]_g  \\
&=[q,(1,1)]_g\xi^g_n((1,\phi^*_{n-1}(q))\Theta_{n-1}(\xi^f_n((u,v))  \\
&=[q,(1,1)]_g\xi^g_n(\varphi^*_{n-1}(q)\Theta_{n-1}([1,(u,v)]_f)  \\
&=[q,(1,1)]_g\ol{\xi^g_n}(q_\nabla((\varphi^*_{n-1}(q))))\Theta_{n-1}([1,(u,v)]_f)  \\
&=[q,(1,1)]_g\ol{\xi^g_n}(q_\nabla((\varphi^*_{n-1}(q))))\Theta_{n-1}([1,(u,v)]_f)  \\
&=[q,(1,1)]_g\ol{\xi^g_n}([1,cg_n(xz\i)]_M)\Theta_{n-1}([1,(u,v)]_f)  \\
&=[q,(1,1)]_g\ol{\xi^g_n}([bp_n(x),cg_n(x)]_M)\Theta_{n-1}([1,(u,v)]_f)  \\
&=[q,(1,1)]_g[1,(bp_n(x),cg_n(x))]_g\Theta_{n-1}([1,(u,v)]_f)  \\
&=\Theta_{n-1}([q,(1,1)]_f)\Theta_{n-1}([1,(u,v)]_f)  \\
&=\Theta_{n-1}([q,(u,v)]_f)
\end{split}
\end{align}

\bibliographystyle{plain}
\bibliography{myrefs}

\end{document}